\documentclass[final,hidelinks,onefignum,onetabnum]{siamart250211}

\usepackage[english]{babel}
\usepackage{color}
\usepackage{amsmath}
\usepackage{amsfonts}
\usepackage{graphicx}
\usepackage{algorithmic}
\usepackage{comment}
\usepackage[normalem]{ulem}
\usepackage{multirow}
\usepackage{booktabs}
\usepackage{indentfirst}
\usepackage{caption}
\usepackage{subcaption}
\usepackage{enumitem}

\headers{Neural Approximate Inverse}{T. Xu, R. Li, and Y. Xi} 

\title{Neural Approximate Inverse Preconditioners}

\author{Tianshi Xu\thanks{Department of Mathematics, University of Kentucky, Lexington, KY 40502, USA (\texttt{\email{tianshi.xu@uky.edu}}). Research of T. Xu is supported by LLNL-LDRD program under Project No. 24-ERD-033.} 
\and Rui Peng Li\thanks{Center  for Applied  Scientific Computing, Lawrence  Livermore National  Laboratory,  P. O.  Box 808,  L-561, Livermore,   CA  94551   {(\texttt{\email{li50@llnl.gov}})}.  This   work  was performed under the  auspices of the U.S. Department  of Energy by Lawrence Livermore National    Laboratory under Contract DE-AC52-07NA27344 (LLNL-JRNL-2012720) and was supported by the LDRD program under Project No. 24-ERD-033 and Project 25-ERD-051.
The United States Government retains, and the publisher, by accepting the article for publication, acknowledges that the United States Government     retains a non-exclusive, paid-up, irrevocable, world-wide license to publish              or reproduce the published form of this manuscript, or allow others to do so, for United States Government purposes.
}\and Yuanzhe Xi
\thanks{Department of Mathematics, Emory University, Atlanta, GA 30322, USA (\texttt{\email{yxi26@emory.edu}}). Research of Y. Xi is supported by NSF DMS-2338904.} 
}
\graphicspath{{./}}
\newcommand{\rev}[1]{{#1}}

\ifpdf\hypersetup{
      pdftitle={Neural Approximate Inverse Preconditioners},
  pdfauthor={T. Xu, R. Li, and Y. Xi}
}
\fi

\begin{document}

\maketitle
\begin{abstract}
In this paper, we propose a learning-based framework for constructing approximate inverse preconditioners for partial differential equations (PDEs) by approximating the Green's function of the underlying differential operator with neural networks (NNs).
The training framework combines four components: an adaptive multiscale neural architecture ($\alpha$MSNN) that resolves hierarchical features across near-, middle-, and far-field regimes; \rev{optional coarse-grid solution anchors that improve identifiability}; a multistage training strategy that progressively refines the Green's function representation across spatial scales; and an overlapping domain decomposition that enables local adaptation while preserving global consistency.
Once trained, the learned Green's function provides a continuous approximation to the inverse operator that can be discretized on a given mesh to construct approximate inverse preconditioners. \rev{Within the domain decomposition framework, selected blocks of the learned inverse can also be used to approximate the inverse Schur complement, avoiding the need to assemble a full dense inverse.}
Numerical experiments on challenging PDE problems demonstrate that the resulting preconditioners \rev{provide robust and scalable convergence and are competitive with classical two-level preconditioning methods}.
\end{abstract}

% REQUIRED
\begin{keywords} 
  Approximate inverse preconditioning, neural networks, domain decomposition, Krylov subspace methods, Schur complement
\end{keywords}

\begin{AMS}
   	65N55,65N80,68T07,65F08,65F10  	
\end{AMS}

% Original submission section: intro.tex
\section{Introduction}\label{sec:intro}
We consider iterative methods for solving  linear system
\begin{equation}
    \mathbf{A} \mathbf{x} = \mathbf{b},
\end{equation}
where the coefficient matrix \(\mathbf{A} \in \mathbb{R}^{n \times n}\) 
is large and sparse,
arising from the discretization of partial differential equations (PDEs), \(\mathbf{x}\) is the discrete solution, and the right-hand-side vector \(\mathbf{b}\) encodes source terms and the boundary conditions.
Krylov subspace methods---such as GMRES, CG, and BiCGSTAB---are widely used for solving such systems. However, the convergence highly depends on the spectral properties of \(\mathbf{A}\). 
When \(\mathbf{A}\) is ill-conditioned, e.g., arising from problems with strong anisotropy, discontinuous or highly variable coefficients,
complex mesh geometries,
or when the operator is indefinite,
these methods may converge slowly or fail to converge altogether \cite{Saad_2003}. 
To mitigate such issues, preconditioning 
techniques are
typically employed to improve spectral properties and accelerate convergence. Applying the left preconditioning approach with the
preconditioner matrix \(\mathbf{M}\) 
yields the transformed system
\[
\mathbf{M}^{-1} \mathbf{A} \mathbf{x} = \mathbf{M}^{-1} \mathbf{b},
\]
in which the preconditioned matrix \(\mathbf{M}^{-1} \mathbf{A}\) ideally exhibits a clustered spectrum and a reduced condition number. Designing preconditioners that are simultaneously effective, scalable, and robust---especially for systems derived from complex or heterogeneous PDEs---remains challenging.

Several classes of preconditioners have been developed to improve the performance of Krylov subspace methods in this context. Among general-purpose strategies, incomplete LU (ILU) factorizations remain widely used due to their algebraic formulation and ease of implementation \cite{chow1997experimental,saad_ilut:_1994, Saad_2003,gpuilu}. Variants incorporating complex shifts \cite{osei-kuffuor_preconditioning_2010,rational} or low-rank corrections \cite{nmllr,mllr,Xu_Kalantzis_Li_Xi_Dillon_Saad_2022} have improved robustness in challenging regimes. However, ILU-based methods require sequential triangular solves, which hinders their parallel scalability. Algebraic multigrid (AMG) methods offer near-optimal complexity for  symmetric positive definite (SPD) systems from 
elliptic PDEs 
\cite{doi:10.1137/S1064827598339402,ruge_algebraic_1987,Xu_Zikatanov_2017}. Nevertheless, the performance of AMG may degrade
with highly heterogeneous coefficients and complex mesh geometries,
and can become ineffective
for indefinite problems.
Moreover, the multilevel property makes it
challenging to identify the optimal AMG algorithm and parameter
setting for a given problem.
\rev{Hierarchical matrix methods exploit low-rank off-diagonal blocks to construct data-sparse representations of matrices and their inverses~\cite{hbook,khoromskij2003datasparse}, with nested-basis $\mathcal{H}^2$ representations providing a further route for kernel compression~\cite{ddh2}.
For large-scale problems, domain decomposition methods divide the problem into smaller subdomain problems while accounting for the coupling between subdomains. Widely used approaches include Schwarz methods \cite{djn2015dd,jolivet2021hpddm} and algebraic Schur complement methods \cite{li_saad_sosonkina_2003_parms,Xu_Kalantzis_Li_Xi_Dillon_Saad_2022}.}

\rev{Machine learning has also been used to approximate PDE solutions and to construct components of iterative solvers. Physics-informed neural networks (PINNs), for example, approximate the solution of a PDE with prescribed coefficients, forcing terms, boundary conditions, and other problem data~\cite{raissi2019physics}. A network trained for one such problem instance therefore represents that particular solution rather than an approximation to the inverse operator, and hence does not directly provide a reusable preconditioner for Krylov methods. Operator-learning approaches, such as DeepONet~\cite{lu2021learning} and the Fourier neural operator (FNO)~\cite{li2020fourier}, instead learn mappings from problem inputs, such as coefficients or forcing terms, to solution fields and can therefore be applied across families of PDE instances. In their supervised formulations, these methods are trained on sampled input--solution pairs, so their performance on new problems depends on how well the training distribution covers the target problem class. Learned solution operators can also be integrated with classical iterative solvers rather than used only as standalone surrogates. DeepONet, for example, has been incorporated into iterative schemes~\cite{zhang2024hints} and used to construct preconditioners~\cite{kopanicakova_karniadakis_2025_deeponet_preconditioning}. In hybrid methods such as HINTS~\cite{zhang2024hints}, learned corrections target error components that are difficult for standard relaxation, while classical iterations such as damped Jacobi reduce complementary components. Other architectures, including convolutional neural networks~\cite{huang2023learning} and graph neural networks~\cite{brarda2026spectral_fsai}, have likewise been used to construct learned smoothers and preconditioners. Although these approaches differ in architecture and in how the learned model is incorporated into the solver, they generally rely on training data sampled from PDEs or matrices intended to represent the target problem distribution.}

In this paper, we consider the family of
approximate inverse preconditioners 
that seek to explicitly approximate the inverse  matrix 
\(\mathbf{A}^{-1}\) \cite{benzi95,Benzi_Tuma_1998,chowpattern,afn}. These preconditioners are well-suited to modern parallel hardware, 
as their applications only require matrix–vector multiplications. The classical sparse approximate inverse techniques typically prescribe the sparsity pattern and use numerical heuristics to compute the preconditioner. 
While being successful in a variety of problems, 
it is known that these methods can fail to capture the nonlocal structure of \(\mathbf{A}^{-1}\), particularly in problems with long-range interactions or highly variable coefficients.

\rev{We exploit a continuous Green's function representation to learn an approximate inverse operator directly from the governing PDE. For a fixed differential operator and boundary conditions, the resulting kernel can be applied to different right-hand sides, including the residual vectors generated during a Krylov iteration. Rather than relying on a large dataset of high-fidelity forcing--solution pairs, we train the kernel using PDE and boundary residuals, with optional coarse-grid solution anchors. The training procedure samples source locations and operator parameters and therefore does not require an application-specific distribution of physical right-hand sides or Krylov residuals. Once trained, the learned Green's function can be discretized on different meshes to produce approximate inverse matrices or selected inverse blocks. These discrete approximations can then be compressed and, when appropriate, further improved through algebraic refinement. The resulting operators can be used as preconditioners for Krylov methods or as approximations to inverse Schur complements, with their quality evaluated through both iterative convergence and the spectral properties of the preconditioned system. Because a fully assembled learned inverse is dense and incurs quadratic storage and matrix--vector multiplication costs, we further consider several structured realizations. These include $\mathcal{H}$-matrix compression, sparse approximate inverse representations, and a two-level ILU framework that uses only the interface block of the learned inverse.}

The main contributions of this work are as follows:
\begin{enumerate}[leftmargin=18pt, label=\textbf{\arabic*}., itemsep=0.5\baselineskip]
\item We develop an ML-based framework for approximating Green's functions with high accuracy and low computational cost. The approach integrates four key components: (i) an adaptive multiscale NN approximation decomposed into near-, middle-, and far-field components, each governed by trainable scaling parameters; (ii) auxiliary anchor data from numerical solutions on smaller grids to improve identifiability; and (iii) a multi-staged training protocol that stabilizes learning and resolves small-scale structure; and (iv) an overlapping domain decomposition approach that captures the locally varying behavior while ensuring a globally consistent representation of the Green's function.
\rev{ \item The learned approximate inverse can be compressed using $\mathcal{H}$-matrix or sparse approximate inverse techniques to reduce storage and matrix-vector multiplication costs. We further incorporate the learned Green's function into a two-level ILU framework to approximate the inverse of the Schur complement. 
Only the values of the learned Green's function at pairs of interface nodes are needed to approximate the inverse of the Schur complement. The matrix entries can be computed independently, allowing parallel construction on distributed-memory systems. The matrix can also be formed on a GPU concurrently with the ILU factorization of the interior blocks on CPUs.
\item We extend our earlier work~\cite{hao2024multiscale} to a broader range of PDEs, including indefinite and highly heterogeneous cases, and explore applications to parametric PDEs. The numerical results demonstrate that our approach significantly accelerates the convergence of Krylov solvers.}
\end{enumerate}

The remainder of the paper is organized as follows. Section~\ref{sec:background} reviews the mathematical background on Green's functions and their relationship to the inverse of discretized differential operators. Section~\ref{sec:method} presents an efficient multiscale framework for learning Green's functions via NNs. \rev{Section~\ref{sec:preconditioner} describes sparse and hierarchical representations and the two-level ILU framework.} Section~\ref{sec:experiments} reports numerical results demonstrating the scalability and spectral effectiveness of the proposed approach. Finally, Section~\ref{sec:conclusion} summarizes the main contributions and discusses directions for future work.

% Original submission section: background.tex
\section{Background}
\label{sec:background}

The connection between the inverse of a discretized coefficient matrix $\mathbf{A}$ and the Green's function of the underlying PDE forms the mathematical foundation of our approach. In this section, we review this relationship and outline related work on data-driven methods for approximating Green's functions.

\subsection{Green's Function}
Consider the linear boundary value problem
\begin{equation}\label{eq:pde}
\begin{cases}
\mathcal{L}u(x) = g_{1}(x), & x\in\Omega,\\[4pt]
\mathcal{B}u(x) = g_{2}(x), & x\in\partial\Omega,
\end{cases}
\end{equation}
where $\mathcal{L}$ denotes a linear differential operator, $\mathcal{B}$ is the boundary operator, $\Omega \subseteq \mathbb{R}^d$ is the computational domain, $u$ is the unknown solution, $g_1$ is a source term, and $g_2$ prescribes the boundary condition.
The associated Green's function $G(x,y)$ satisfies
\begin{equation}\label{eq:green}
\begin{cases}
\mathcal{L}_{x} G(x,y) = \delta(x-y), & x,y \in \Omega,\\[4pt]
\mathcal{B}_{x} G(x,y) = 0, & x \in \partial\Omega,
\end{cases}
\end{equation}
 where the subscript $x$ indicates that the operators act on the $x$-variable of $G(x,y)$,
and
$\delta(x-y)$ is the Dirac delta distribution, characterized by 
\begin{equation}
\delta(x) = 0, \ \forall x \neq 0 \quad \mbox{and} \quad 
\int_{\mathbb{R}^d} \delta(x)\,dx = 1.
\end{equation}
Once the Green's function is known, the solution of \eqref{eq:pde} can be expressed as an integral representation. 
For instance, consider a second-order elliptic PDE with Dirichlet boundary condition,
\begin{equation}\label{eq:model}
-\nabla \!\cdot\! \bigl( a(x)\nabla u(x) \bigr) = g_{1}(x), \qquad
u|_{\partial\Omega} = g_{2}(x),
\end{equation}
where $a(x)$ is a positive scalar diffusion coefficient.
Green's function $G(x,y)$ satisfies \emph{homogeneous} boundary conditions and
the solution admits the representation
\begin{equation}\label{eq:integral_PDE}
u(x) = \int_{\Omega} G(x,y)\, g_1(y)\, d y
- \int_{\partial\Omega} g_{2}(y)\, a(y)\,
\bigl(\nabla_{y} G(x,y) \!\cdot\! n_{y}\bigr)\, dS(y).
\end{equation}

\subsection{Connection to Matrix Inverses}
\begin{figure}
    \centering
    \includegraphics[width=0.95\linewidth]{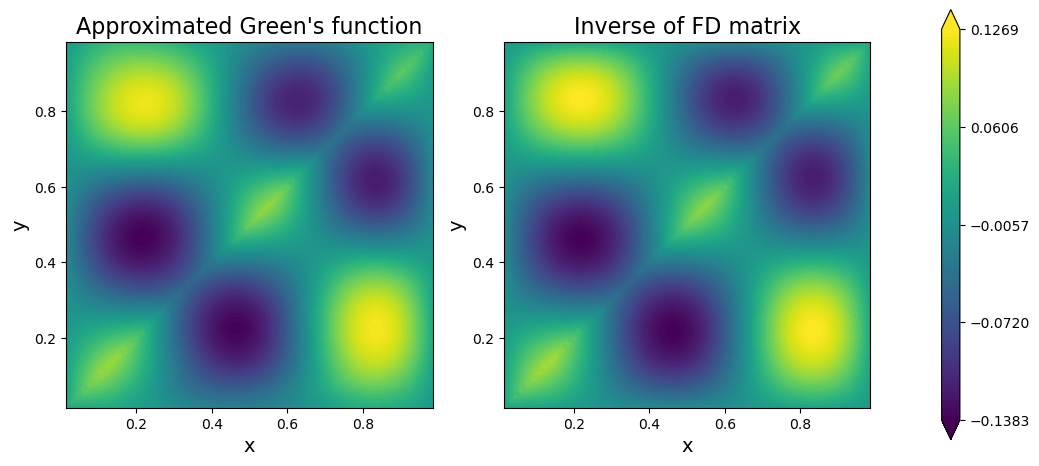}
    \caption{Comparison between an approximation of the Green's function and the inverse of the stiffness matrix.
    \textit{Left:} A NN approximation of the Green's function of the operator $\mathcal{L}=-\Delta - 50$ obtained, evaluated on a uniform grid over the domain $[0,1]$ with $64$ interior points.
    \textit{Right:} Inverse of the finite difference discretization of the same operator on the same grid.}
    \label{fig:inv_vs_green}
\end{figure}

As indicated in \eqref{eq:green}, 
the Green's function serves as the kernel of the continuous inverse operator. 
Upon
discretization, 
the operator 
$\mathcal{L}$ is replaced by the stiffness matrix $\mathbf{A}$
and its inverse $\mathbf{A^{-1}}$ plays the same role as $G(x, y)$.
Hence, the discrete Green's function can be viewed as an
approximation to $\mathbf{A}^{-1}$. 
To make this correspondence precise,
consider a finite difference discretization with interior grid nodes $\{x_i\}_{i=1}^{n} \subset \Omega$ and boundary nodes $\{x_k\}_{k=n+1}^{n+m} \subset \partial\Omega$, 
with the mesh step size $h$. 
Here we multiply the pointwise finite-difference equations by $h^d$, so that $\mathbf A=h^d\mathbf L_h$, where $\mathbf L_h$ is the interior finite-difference operator.
Discretizing PDE \eqref{eq:model} yields the block linear system
\begin{equation}
\label{eq:AB}
\begin{bmatrix}
\mathbf{A} & \mathbf{B} \\
\mathbf{0} & \mathbf{I}
\end{bmatrix}
\begin{bmatrix}
\mathbf{x} \\
\mathbf{x}_{\Gamma}
\end{bmatrix}
=
\begin{bmatrix}
\mathbf{f} \\
\mathbf{g}_{2}
\end{bmatrix},
\quad
f_i = h^d g_1(x_i),\quad i=1,\dots,n,
\end{equation}
where $\mathbf{A}\in\mathbb{R}^{n\times n}$ couples  the interior nodes, $\mathbf{B}\in\mathbb{R}^{n\times m}$ couples the interior and 
boundary, and $\mathbf{g}_2$ encodes the boundary condition. Eliminating $\mathbf{x}_{\Gamma}$ 
from \eqref{eq:AB}
leads to reduced interior system
\begin{equation}
\mathbf{A} \mathbf{x} = \mathbf{f} - \mathbf{B}\mathbf{g}_2.
\label{eq:interiorSystem}
\end{equation}
Analogously, discretizing the Green's function 
representation \eqref{eq:integral_PDE} by midpoint quadrature gives, 
for an interior node $x_i$, we have
\begin{equation}
\hat u(x_i) = \sum_{j=1}^{n} G(x_i, x_j)\,g_1(x_j)\,h^d
- \sum_{k=1}^{m} \frac{\partial G}{\partial n_y}(x_i, x_{n+k})\,a(x_{n+k})\,g_2(x_{n+k})\,w^{\Gamma}_k,
\label{eq:discInt}
\end{equation}
where $w_k^\Gamma \sim h^{d-1}$ are boundary quadrature weights. 
Introducing the matrices
\begin{equation}
\mathbf{G}_{ij} = G(x_i, x_j), \quad
\mathbf{D}_{ik} = \frac{\partial G}{\partial n_y}(x_i, x_{n+k}),
\end{equation}
and the diagonal matrix
$\mathbf{W}_{\Gamma} = \operatorname{diag}(a(x_{n+1})w^{\Gamma}_1, \dots, a(x_{n+m})w^{\Gamma}_m)$, 
we can write
\begin{equation}
\mathbf{\hat u} = (h^d\mathbf{G})\,\mathbf{g}_1 - (\mathbf{D}\mathbf{W}_{\!\Gamma})\,\mathbf{g}_2.
\label{eq:discRep}
\end{equation}
Since both $\mathbf{\hat u}$ and $\mathbf{x}$ approximate 
$\mathbf{u}$, 
we obtain the approximate relation
\begin{equation}\label{eq:match}
(h^d\mathbf{G})\mathbf{g}_1 - (\mathbf{D}\mathbf{W}_{\!\Gamma})\mathbf{g}_2
\;\approx\; \mathbf{A}^{-1}\!\left(h^d \mathbf{g}_1 - \mathbf{B}\mathbf{g}_2\right),
\end{equation}
from which comparing the coefficients of $\mathbf{g}_1$, it follows that
\begin{equation}
\mathbf{G} \;\approx\; \mathbf{A}^{-1}.
\end{equation}
Hence, the discrete Green's function provides an approximation to the inverse 
of the coefficient matrix. 
Figure~\ref{fig:inv_vs_green} compares an approximation of the Green's function (left) with the finite-difference inverse $\mathbf A^{-1}$ (right) for $\mathcal L=-\Delta-50$ on the same uniform grid on $[0,1]$.
The comparison of these two figures confirms that the Green's function closely resembles the inverse matrix. A similar correspondence arises in finite element discretizations, where the Green's operator relates to the inverse of the stiffness matrix combined with mass and trace operators. 

\subsection[Related work]{Related work}\label{sec:related-work}

Although Green's functions provide a natural representation of inverse operators, they are rarely available in closed form except for simple domains and operators.
\rev{DeepGreen~\cite{gin2021deepgreen} learns a transformed Green representation from forcing--solution pairs, while BI-GreenNet~\cite{lin2023bi} uses a boundary integral representation. GreenONet of Aldirany et al.~\cite{aldirany2024greenonet} combines a Green integral representation with physics-informed training. The Green's operator network (GON) of Gu et al.~\cite{gu2025greensonet} learns Green's functions and their gradients to compute solutions through volume and boundary integrals. DecGreenNet~\cite{wimalawarne2025decgreennet} learns a low-rank approximation of the Green's function to reduce the cost of numerical integration. Green Multigrid Network~\cite{lin2024green} combines a piecewise neural representation with multilevel integration. For parametric Green operators, Feliu-Fabà et al.~\cite{feliu_faba2020meta_pdo} learn compressed wavelet representations through meta-learning. Teng et al.~\cite{teng2022learning} partition the source domain and train a separate physics-informed network for each source subdomain.

Learned Green representations have also been used directly in iterative solution, including earlier discrete Green neural methods~\cite{ichimura2020green_nn,tang2022neural_green}. Sun et al.~\cite{sun2026singularity_green} encode Green singularities and deploy learned matrices as dense preconditioners, two-level overlapping Schwarz coarse corrections, and components of hybrid iterations. Neural Green's Operators (NGOs) learn finite-dimensional Green coefficient representations and support preconditioner construction~\cite{melchers2024neural}; a follow-up study uses NGOs to construct a fixed-basis Galerkin coarse correction~\cite{melchers2026coarse_solve}. These methods provide explicit routes from learned Green information to preconditioners and iterative corrections. 

Our approach focuses on adaptive multiscale learning of Green's functions and their use in two-level ILU preconditioning, where the learned Green's function on the interface provides an approximation to the inverse Schur complement.\par}

% Original submission section: greenlearning.tex
\section{Efficient Green's Function Learning}\label{sec:method}

The effectiveness of using NNs to approximate the Green's function
and thereby construct a preconditioner relies critically 
on the accuracy of representing the true inverse operator. 
Achieving the accuracy is nontrivial: it requires the careful co-design of the NN architecture, loss function, and training strategy. If these components are not aligned, the training may converge to a suboptimal solution that fails to capture essential features of the underlying PDE---such as anisotropy, boundary layers, or localized singularities. Such inaccuracies undermine the spectral effectiveness of the resulting preconditioner and ultimately impair solver performance.

\begin{figure}[htbp]
\centering
\includegraphics[width=0.95\textwidth]{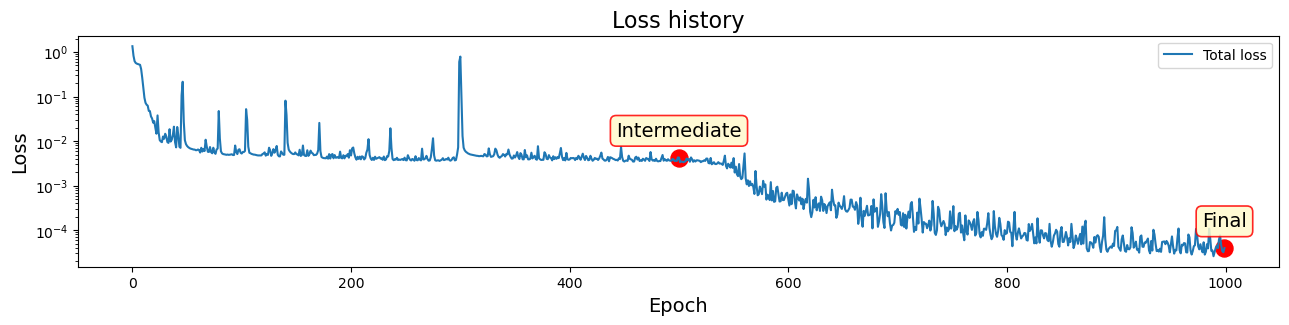}
\includegraphics[width=0.95\textwidth]{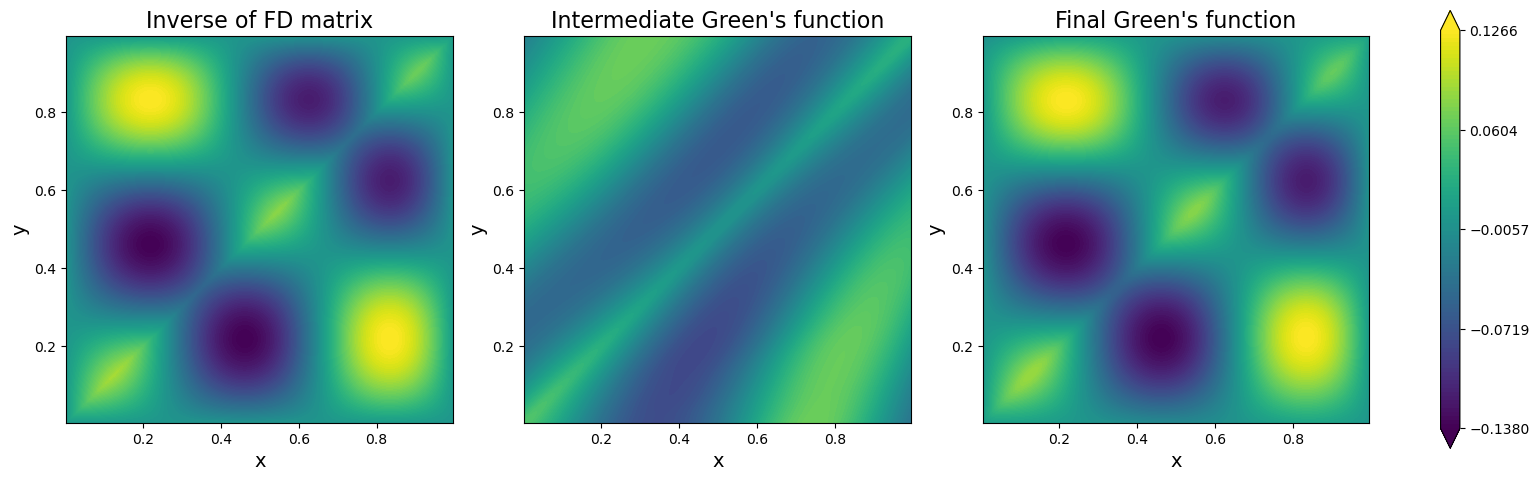}
\caption{Effect of approximation accuracy on the prediction accuracy of the learned Green's function for the convection–diffusion-reaction problem~\eqref{eq:conv-diff}. 
Approximations are compared against a numerical reference solution. 
\textit{Top:} training loss history.
\textit{Bottom left:} numerical reference solution.
\textit{Bottom center:} training to a relative loss of $10^{-2}$ yields large errors. 
\textit{Bottom right:} training to $10^{-4}$ substantially reduces the error and recovers the expected reaction-induced waves.}
\label{fig:motivation_accuracy}
\end{figure}

To illustrate the impact of approximation quality, we consider the following convection–diffusion-reaction equation defined on the domain $[0,1]$:
\begin{equation}\label{eq:conv-diff}
\mathcal{L}u = -\nabla \cdot (a(x)\nabla u) + {b}(x) \cdot \nabla u + c(x)u.
\end{equation}
{Figure~\ref{fig:motivation_accuracy} illustrates the qualitative evolution of the learned Green's function throughout the training process. 
The top panel shows the training loss trajectory of a 3-layer multilayer perceptron (MLP), 
with two selected checkpoints highlighted: an intermediate model where the loss plateaus 
around $10^{-2}$, and the final model trained to a tolerance of $10^{-4}$. 
The bottom row compares the ground truth inverse operator (left) with the network output at these two stages. 
At the intermediate checkpoint (middle panel), the learned function remains 
qualitatively incorrect, failing to capture the key structural features of the true operator. 
In contrast, after further training to a much lower loss, 
the final model (right panel) successfully reproduces the characteristic 
anisotropic patterns of the reference solution. 
This progression demonstrates that 
achieving high numerical accuracy in approximating Green's functions
is essential for preserving the correct physical properties, 
a vital requirement when such operators are used for preconditioning.
\rev{Establishing a quantitative relationship between the learned Green-kernel error and the spectrum of the preconditioned discrete system remains an important direction for future work. In particular, the sampled training loss does not by itself control the error of the resulting discrete approximate inverse, which also depends on the discretization and properties of the underlying PDE. In practice, preconditioning quality can be assessed using a small set of independent validation problems and right-hand sides, with true relative residuals and Krylov iteration counts providing direct measures of solver performance.}}

However, accurately learning Green's functions for PDEs---especially those with anisotropic coefficients or multiscale structure---remains a significant challenge. Standard NN architectures often fail to resolve the sharp transitions and singularities.
We next describe the key components of the proposed framework,
which include an adaptive multiscale architecture, a strategy for selecting auxiliary anchor data, and a multi-$\varepsilon$ staged training protocol coupled with an overlapping domain decomposition method that can capture
the local variations while preserving the global consistency in the Green's function.

\subsection{Multiscale Neural Networks}
Accurate approximations to Green's functions require resolving behavior across a wide range of spatial scales. Near the source point \(y\) 
(i.e., $x\approx y$),
the kernel \(G(x, y)\) exhibits singular and sharply peaked features. 
In contrast, in the far field (i.e., for large $\lvert x-y \rvert$), 
it decays smoothly and admits a low-rank structure. 
The contrast between low regularity near the source and smoother behavior away from it can make single-scale neural approximations inefficient~\cite{hao2024multiscale}.
To overcome this limitation, we employ 
a multiscale neural architecture specifically designed to reflect the hierarchical structure of such target functions.

In practice, the learning task can be made more tractable by first regularizing the Dirac delta in \eqref{eq:green} using a Gaussian kernel with
a small width $\varepsilon > 0$:
\begin{equation}
\delta(x-y) \approx \mathcal{N}_\varepsilon(x,y)=\left(\frac{1}{\varepsilon\sqrt{\pi}}\right)^d\exp\left(-\frac{\|x-y\|_2^2}{\varepsilon^2}\right).
\label{eq:delta}
\end{equation}
This transforms the original problem into one with a regularized Green's function \(G_\varepsilon\):
\begin{equation}
\begin{cases}
    \mathcal{L}_x G_\varepsilon(x,y) = \mathcal{N}_\varepsilon(x,y), & x,y \in \Omega,\\
    \mathcal{B}_x G_\varepsilon(x,y) = 0, & x \in \partial\Omega, y \in \Omega.
\end{cases}
\label{eq:regularized_green}
\end{equation}
Although this regularization removes the singularity at \(x = y\), the Gaussian source term \(\mathcal{N}_\varepsilon(x,y)\) remains sharply peaked and poorly regular as \(\varepsilon \to 0\). Consequently, the resulting Green's function \(G_\varepsilon(x,y)\) still exhibits strong spatial variation near the source and long-range decay away from it. Capturing this behavior demands a network architecture that can resolve localized structure without sacrificing global accuracy.
The Multiscale Neural Network (MSNN) framework proposed in \cite{hao2024multiscale} addresses the challenge of learning Green's functions by decomposing the approximation into two 
components with distinct spatial resolutions: the low-level 
component targets at  localized, high-amplitude features near the 
source, while the high-level one captures 
the smoother and global structure.
The NN is represented in the general form
\begin{equation} \label{eq:general_msnn}
G_{\text{MSNN}}(x, y; \theta) = \phi_1(x, y; \theta_1) + \phi_2(x, y; \theta_2),
\end{equation}
with each component  \(\phi_i\), $i=1,2$, defined as a \emph{scaled} feed-forward NN
\begin{equation} \label{eq:phii}
\phi_i(x, y; \theta_i) = \varepsilon^{\alpha_i} \sum_{j=1}^{n_i} a_{ij} \, \sigma\left( \frac{w_{ij} \cdot (x, y) + b_{ij}}{\varepsilon^{\beta_i}} \right),
\end{equation}
where $\sigma$ denotes the nonlinear activation function, 
$n_i$ is the number of neurons, and 
$w_{ij}$ and $b_{ij}$ are the weight and bias terms, respectively. 
The scaling parameters $\alpha_i$ and $\beta_i$ 
control the amplitude and spatial frequency of the resulting basis functions.
MSNN provides a principled framework for modeling Green's functions in Poisson-like problems, capturing both the singular near-field behavior around the source and the smooth far-field structure. 
However, the architecture faces two notable limitations. First, it is restricted to just two branches---a near-field and a far-field component---making it unable to represent the middle scales that often arise between sharp source regions and the global background in more general PDEs. Second, each branch relies on fixed scaling parameters $(\alpha_i,\beta_i)$. While sufficient for simple, homogeneous settings, these fixed scales often fail to adapt to spatial variations in heterogeneous media, thereby reducing the robustness and generality of the method.

\subsection{Adaptive Multiscale Neural Networks}

To address these limitations, we propose an \emph{Adaptive Multiscale NN} ($\alpha$MSNN) framework. The architecture retains the basic near- and far-field decomposition, but augments it with additional middle branches and replaces the fixed scaling parameters $(\alpha,\beta)$ with trainable, source-dependent exponents. This added flexibility enables the network to capture residual structures across multiple scales.

To motivate the design of $\alpha$MSNN, we analyze a model setting 
where the local scaling of the Green's function can be made explicit. Consider the operator $\mathcal{L}$ in \eqref{eq:conv-diff} with a fixed source $y_0$. 
We seek the behavior of $G(x, y_0)$ when
$x$ lies
in an distance-$\varepsilon$ neighborhood of $y_0$, i.e., $x = y_0 + \varepsilon \xi$.
For sufficiently small $\varepsilon$, the coefficients admit 
$a(x)\approx a(y_0)$, $b(x)\approx b(y_0)$, and $c(x)\approx c(y_0)$.
Under the change of variables, the derivatives rescale as
$\nabla = \varepsilon^{-1}\nabla_\xi$
and
$\Delta = \varepsilon^{-2}\Delta_\xi$.
Using the product rule
\begin{equation}
-\nabla\!\cdot\!\big(a(x)\nabla \big)
= -a(x)\,\Delta - \nabla a(x) \!\cdot\!\nabla,
\end{equation}
and freezing $a$ in the second-order term while retaining the first-order in the
convection term, we obtain the locally rescaled operator
\begin{align}
\mathcal{L}_x
&\approx -\varepsilon^{-2} a(y_0)\nabla_\xi\!\cdot\!\nabla_\xi
+
\varepsilon^{-1}\big(b(y_0)-\nabla a(y_0)\big)\!\cdot\!\nabla_\xi
+
c(y_0).
\label{eq:scalings}
\end{align}
Thus the diffusion, {convection}, and reaction contributions appear at distinct orders in $\varepsilon$.
Consequently, any accurate local approximation of the Green's function must capture this multi-scale structure.

\subsubsection{Near-field Parameterization}
To model the locally singular behavior of $G(x, y_0)$, 
we use the relative displacement $x-y_0$ to express the variation around the source and introduce the spatial scaling below.
To achieve this, we re-center the first-layer parameters of the NN so that its activation depends on $x - y_0$ rather than $x$ and $y_0$ separately. 
The resulting near-field branch is parameterized as
\begin{equation} \label{eq:near-field}
\phi_1(x, y_0; \theta_1)
= \varepsilon^{\alpha_1} \sum_{j=1}^{n_1} a_{1j} \,
\sigma\!\left( \frac{\tilde{w}_{1j} \cdot (x - y_0) + \tilde{b}_{1j}}{\varepsilon^{\beta_1}} \right)
:= \varepsilon^{\alpha_1} f\!\left(\frac{x - y_0}{\varepsilon^{\beta_1}}\right),
\end{equation}
where $\beta_1$ controls the effective spatial width of the activation
and $\alpha_1$ sets its amplitude.
To express the dependence explicitly on $x - y_0$,
parameters are re-centered by decomposing each weight
$w_{1j} = (w_{1j}^x, w_{1j}^y)$ and defining 
\begin{equation}
\tilde{w}_{1j} = w_{1j}^x,
\qquad
\tilde{b}_{1j} = (w_{1j}^x + w_{1j}^y)^\top y_0 + b_{1j}.
\end{equation}
This reparameterization ensures that each activation in $\phi_1$ depends only on the scaled displacement 
$(x - y_0)/\varepsilon^{\beta_1}$, 
so that the branch remains localized around the source location~$y_0$.
In practice, $\sigma$ is implemented using a smooth, saturating activation function (e.g., $\tanh$), 
whose derivative is concentrated within a \emph{narrow transition region} and decays exponentially outside it due to the boundedness of the weights and biases~\cite{hao2024multiscale}. 
Within this region, the network produces a sharply localized response, 
confined to an $\mathcal{O}(\varepsilon^{\beta_1})$ neighborhood of $y_0$, 
thereby aligning the branch with the near-field of the Green's function and suppressing spurious long-range effects. 
This construction also ensures stability under differentiation by $\mathcal{L}$, 
unlike non-decaying activations (e.g., ReLU) or polynomial bases that can propagate long-range artifacts.

To determine suitable values for $\alpha_1$ and $\beta_1$, 
we analyze the asymptotic scaling of the near-field component $\phi_1$ in the neighborhood of $y_0$. 
Introduce the orders (in $\varepsilon$) of the local magnitudes
\begin{equation}
|a(y_0)|\sim \varepsilon^{p_a},\qquad
|b(y_0)-\nabla a(y_0)|\sim \varepsilon^{p_b},\qquad
|c(y_0)|\sim \varepsilon^{p_c}.
\end{equation}
With $\phi_1$ as in \eqref{eq:near-field}, the chain rule gives
\begin{equation}
\nabla \phi_1=\varepsilon^{\alpha_1-\beta_1}(\nabla f)\!\left(\tfrac{x-y_0}{\varepsilon^{\beta_1}}\right),
\qquad
\Delta \phi_1=\varepsilon^{\alpha_1-2\beta_1}(\Delta f)\!\left(\tfrac{x-y_0}{\varepsilon^{\beta_1}}\right).
\end{equation}
Substituting these into $\mathcal{L}_x \phi_1$ yields a refined asymptotic estimate of the terms in~\eqref{eq:scalings}:
\begin{equation}\label{eq:Lphi1-scaling}
\mathcal{L}_x \phi_1 \sim
\begin{cases}
\varepsilon^{p_a + \alpha_1 - 2\beta_1}, & \text{if diffusion dominates},\\[3pt]
\varepsilon^{p_b + \alpha_1 - \beta_1},  & \text{if effective convection dominates},\\[3pt]
\varepsilon^{p_c + \alpha_1},            & \text{if reaction dominates}.
\end{cases}
\end{equation}
For the regularized Green's equation, the right-hand side Gaussian function $\mathcal{N}_\varepsilon$ scales as $\varepsilon^{r}$.
When $\varepsilon\ll 1/\sqrt{\pi}$, $\mathcal{N}_\varepsilon\sim \varepsilon^{-d}$.
To ensure consistency at the leading order, the dominant term in~\eqref{eq:Lphi1-scaling} 
must balance the source, 
which yields
\begin{equation}
\alpha_1=
\begin{cases}
r - p_a + 2\beta_1, & \text{if diffusion dominates},\\[3pt]
r - p_b + \beta_1,  & \text{if effective convection dominates},\\[3pt]
r - p_c,            & \text{if reaction dominates}.
\end{cases}
\end{equation}
Moreover, since $\mathcal{N}_\varepsilon$ has effective width $\mathcal{O}(\varepsilon)$, it is natural to set $\beta_1=1$.

\subsubsection{Far- and middle-field branches}
\begin{figure}[htbp]
    \centering
    \includegraphics[width=0.95\linewidth]{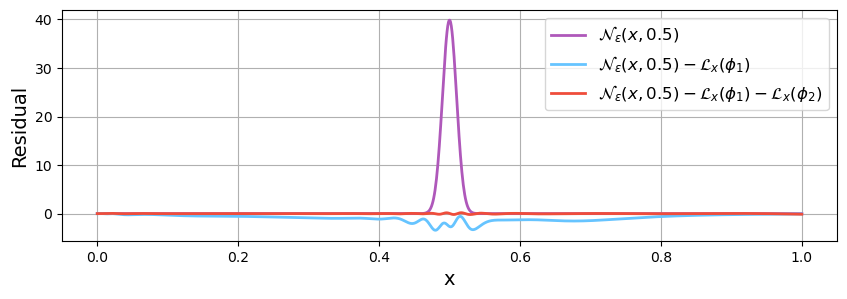}
    \caption{Gaussian source $\mathcal N_\varepsilon$, residual $\mathcal N_\varepsilon-\mathcal L_x\phi_1$, and residual $\mathcal N_\varepsilon-\mathcal L_x(\phi_1+\phi_2)$ for the operator $\mathcal{L}u=-\Delta u$ on $[0,1]$. The progression illustrates how successive components cancel localized structure and refine the multiscale hierarchy at a fixed $y_0=0.5$.}
    \label{fig:gauss residual}
\end{figure}

Away from the source, elliptic regularity guarantees that the Green's function is smooth. This motivates a globally supported 
\emph{far-field} branch,
indexed by $L$,
which follows the same form as \eqref{eq:near-field} but with
$(\alpha_L,\beta_L)=(0,0)$. 
Between this smooth background and the sharply localized near-field response, 
we introduce \emph{middle-field} branches that adaptively capture residual structure 
at intermediate distances, typically
around $|x-y_0|\sim \varepsilon^{\beta_j}$, $0 < \beta_j < 1$.
Each branch employs a smooth, saturating activation (such as $\tanh$),
whose transition region defines its effective spatial scale.
By iteratively applying these branches across multiple scales,
the network constructs a hierarchy in which residuals are progressively canceled,
yielding an increasingly accurate multiscale representation of the Green's function
(see Fig.~\ref{fig:gauss residual}).

The scaling parameters $(\alpha_j,\beta_j)$ for these branches cannot be prescribed \emph{a priori} and must instead be learned from data. 
Once $\phi_1$ is fitted, the residual near $|x-y_0|\sim\varepsilon^{\beta_1}$ generally contains contributions from multiple scales, making fixed choices for subsequent branches unreliable. 
We impose mild ordering constraints: $\beta_j$ decreases monotonically while $\alpha_j$ increases. For $0<\varepsilon<1$, this gives successive branches broader input scales and smaller amplitude prefactors.
Together, these mechanisms provide the $\alpha$MSNN with a systematic strategy to capture features ranging from the singular near-field core to the smooth far-field background.

\subsubsection{Adaptive Multiscale Architecture}
The preceding analysis focused on a fixed source  $y=y_0$ 
to make the scaling behavior transparent. 
In practice, however, the Green's function must be approximated as a bivariate kernel $G(x,y)$ defined for all source points $y$. 
To achieve this, 
$\alpha$MSNN applies the same compositional structure with activations of the form
\begin{equation}
\phi_j(x,y) = \varepsilon^{\alpha_j(y)} \hat{f}\!\left(\frac{x-y}{\varepsilon^{\beta_j(y)}},y\right),
\end{equation}
{with $\hat{f}(\cdot,y_0):=f(\cdot)$.}
{Despite its appeal in diffusion-dominated Poisson settings, hand-picked exponent pairs $(\alpha_j(y),\beta_j(y))$ are often unreliable in practice.
In heterogeneous media, these quantities can vary by orders of magnitude across $y$, causing the dominant physical mechanism to shift throughout the domain.
Moreover, the intermediate (“middle-field”) behavior is difficult to predict.
For these reasons, we let both the amplitude exponent $\alpha_j(y)$ and the spatial exponent $\beta_j(y)$ be learnable functions of the source coordinate $y$.
}
The use of smooth saturating activations such as $\tanh$ ensures that each branch remains localized to the 
$\mathcal{O}(\varepsilon^{\beta_j(y)})$ neighborhood of the source 
while decaying rapidly outside. 
In this way, the network generalizes the fixed-source analysis at $y=y_0$ to a full bivariable representation $G(x,y)$, 
with scale parameters that adapt to spatial variations in the underlying operator.

Finally, 
we can define the complete multiscale $\alpha$MSNN 
approximation with $L$ hierarchical branches as 
\begin{equation}
G_{\alpha\mathrm{MSNN}}(x, y; \theta)  = \sum_{j=1}^L \varepsilon^{\alpha_j(y)} 
\, \hat{f}_j\!\left(\frac{x-y}{\varepsilon^{\beta_j(y)}},y;\,\theta_j\right),
\end{equation}
where each $\hat{f}_j$ represents the functional form of the $j$th branch,
implemented using smooth, saturating activations (e.g., $\tanh$ or its variants)
and parameterized by branch-specific weights~$\theta_j$.
For the near- and middle-field branches, $\beta_j$ decreases with increasing~$j$ while $\alpha_j$ increases, giving progressively broader input scales and smaller amplitude prefactors.
Together, these components form a flexible, data-driven multiscale kernel
that faithfully captures the behavior of Green's functions---from their singular
near-field core to their smooth far-field background.

\begin{figure}[htbp]
    \centering
    \includegraphics[width=0.95\linewidth]{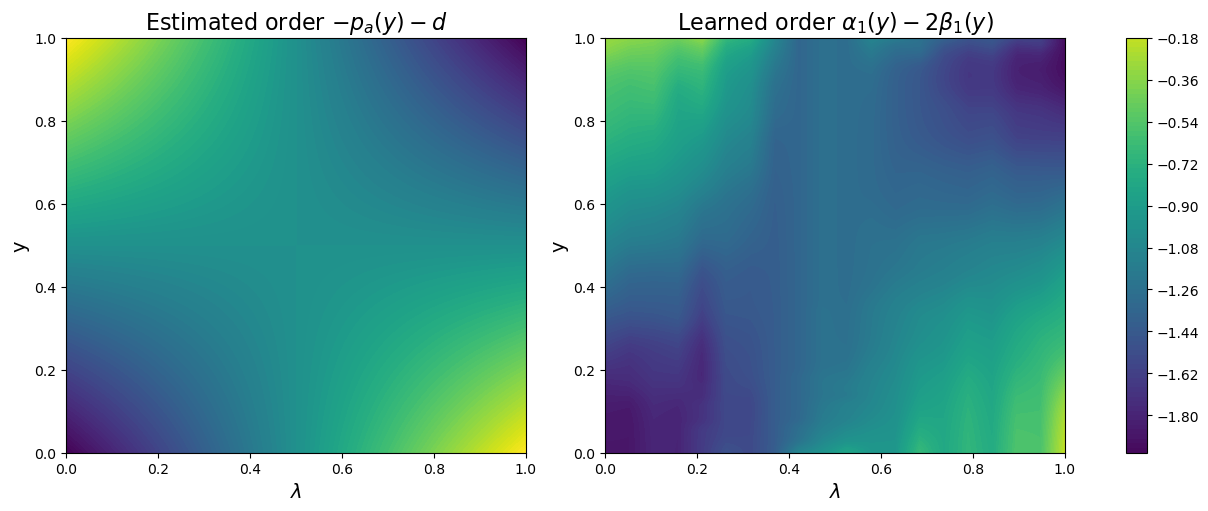}
    \caption{{Learned scaling exponents $(\alpha_1,\beta_1)$ for the near-field branch $\phi_1$ in the 1D diffusion problem with $a(x)=\varepsilon^{(2x-1)(2\lambda-1)}$ and $b(x)=c(x)=0$ with $\varepsilon=0.05$.
For each $\lambda$, the network learns $\alpha_1(y)$ and $\beta_1(y)$ such that $\alpha_1(y)-2\beta_1(y)$ matches the structure of 
the theoretical scaling $-p_a(y)-d$, confirming that $\phi_1$ captures the correct near-field behavior and adapts to spatially varying diffusion.}}
\label{fig:scaling_balance}
\end{figure}

Figure \ref{fig:scaling_balance} illustrates the scaling exponents learned by the near-field branch $\phi_1$ for the one-dimensional diffusion-dominated problem, with the diffusion coefficient parameterized as $a(x)=\varepsilon^{(2x-1)(2\lambda-1)}$ and $b(x)=c(x)=0$.
The parameter $\lambda$ controls the spatial variation of $a(x)$, allowing us to test how the learned exponents respond to different diffusion profiles.
For each fixed $\lambda$, the NN learns $\alpha_1(y)$ and $\beta_1(y)$ from data.
The leading contribution in $\mathcal{L}_x[\phi_1]$ scales as $\varepsilon^{p_a+\alpha_1-2\beta_1}$, implying that a correctly learned branch should satisfy the balance relation $p_a(y)+\alpha_1(y)-2\beta_1(y)\approx -d$. The left panel shows the theoretical diffusion scaling $-p_a(y)-d$, while the right panel shows the corresponding learned quantity $\alpha_1(y)-2\beta_1(y)$ for various $\lambda$.
Their close agreement across $\lambda$ and $y$ demonstrates that the network has successfully captured the near-field scaling behavior and adapted to the local diffusion strength.
The variation of both quantities with $\lambda$ and $y$ also highlights why prescribing fixed exponents is insufficient in heterogeneous media, and confirms that the learned $(\alpha_1,\beta_1)$ provide a data-driven and spatially adaptive representation of the singularity Green's functions.

\subsection{Auxiliary Data from Coarse-grid Solutions}

\leavevmode\rev{Residual-based training alone may not always provide sufficient information to accurately recover the Green's function. To improve the robustness of the training process, we optionally augment the physics-based losses with a small number of approximate Green's function values obtained from a coarse-grid discretization. These values serve as auxiliary anchors during offline training and are not required when constructing or applying the trained preconditioner.}

To see why such additional information can be useful, consider an ill-conditioned differential operator. Although $\alpha$MSNN is sufficiently expressive to approximate the Green's function, residual minimization may be insensitive to certain error components. As a simple illustration, consider the shifted Laplacian 
\begin{equation}
\mathcal{L}u = -u_{xx} - \kappa^2 u, \quad \kappa^2 > 0,
\end{equation}
on $(0,1)$ with homogeneous Dirichlet boundary conditions. \rev{Assume that $\kappa^2$ is not a Dirichlet eigenvalue, so that $\mathcal{L}$ is invertible.} When $\kappa^2$ is close to $\lambda_m=(m\pi)^2$, the corresponding mode $v_m(x)=\sin(m\pi x)$ is a \rev{near-null} mode, since
\begin{equation}
% Revision color removed.
\mathcal{L}v_m=(\lambda_m-\kappa^2)v_m.
\end{equation}
A kernel error of the form $e(x,y)=v_m(x)c(y)$ therefore produces a residual scaled by the small eigenvalue gap $\lambda_m-\kappa^2$. Residual minimization can be insensitive to such error components, making accurate approximation difficult within a finite training budget.

\begin{figure}[!htbp]
    \centering
    \includegraphics[width=0.95\linewidth]{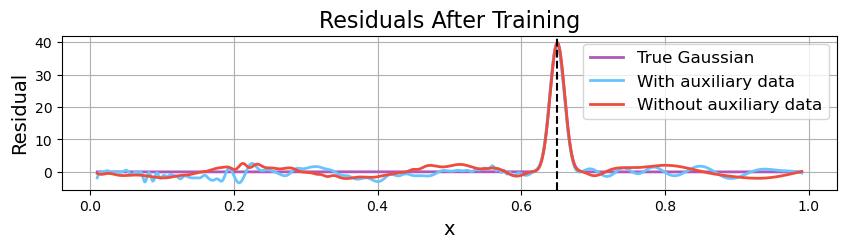}
    \includegraphics[width=0.95\linewidth]{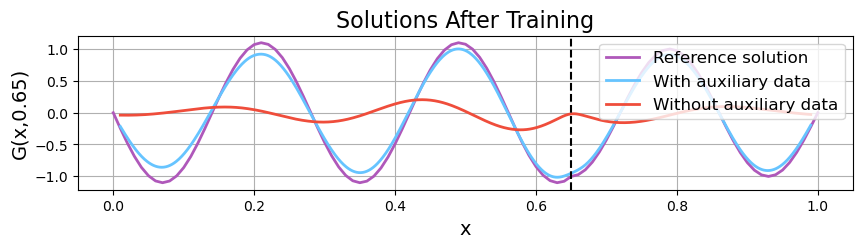}
    \caption{Role of anchor data for $\mathcal L=-0.1\Delta-50$ on $[0,1]$. Both models use two branches with four hidden layers of $50$ neurons each and are trained for $5000$ Adam iterations at $y=0.65$, with Gaussian source standard deviation $0.01$. \textit{Top:} Differential-operator outputs and the Gaussian source. \textit{Bottom:} Learned Green's function slices and the coarse reference.}
    \label{fig:full_comparison}
\end{figure}

In principle, the training problem departs from the ideal continuous formulation in two ways.  
First, the Green's function is approximated by  $G_{\text{$\alpha$MSNN}}(x,y;\theta)$ rather than the exact kernel $G(x,y)$. 
Second, the singular source $\delta(x-y)$ is replaced by a narrow Gaussian $\mathcal{N}_\varepsilon(x,y)$, yielding the continuous residual functional
\begin{equation}
\mathcal{R}_{\text{PDE}}[G_{\alpha\mathrm{MSNN}}] 
= \bigl\|\mathcal{L}_x G_{\alpha\mathrm{MSNN}}(x,y;\theta) - \mathcal{N}_\varepsilon(x,y)\bigr\|_{L^2(\Omega\times\Omega)}^2.
\end{equation}
During training, however, this loss is evaluated only on a finite set of collocation points $S_{\text{res}} \subset \Omega\times\Omega$,
giving the
empirical residual loss
\begin{equation}
\mathcal{L}_{\text{res}}[G_{\alpha\mathrm{MSNN}}] 
= \frac{1}{|S_{\text{res}}|}
\sum_{(x,y)\in S_{\text{res}}} 
\Bigl|\mathcal{L}_x G_{\alpha\mathrm{MSNN}}(x,y;\theta) - \mathcal{N}_\varepsilon(x,y)\Bigr|^2.
\end{equation}
Finite collocation and optimization errors can further limit the accuracy, so a small sampled residual need not imply an accurate Green's function (see Fig.~\ref{fig:full_comparison}). To provide optional additional supervision, we introduce a small set of \emph{anchor points}
\begin{equation}
\{(x_i^{\!*}, y_i^{\!*}, g_i^{\!*})\}_{i=1}^m,
\end{equation}
where $g_i^{\!*}$ are approximate samples of the true Green's function obtained from a coarse finite-difference discretization. 
The connection between Green's functions and matrix inverses described in Section~2.2 provides a simple way to obtain these values from the inverse of a coarse-grid discretization matrix.
The auxiliary (anchor) loss is then defined as
\begin{equation}
\mathcal{L}_{\text{aux}}[G_{\alpha\mathrm{MSNN}}] 
= \frac{1}{m}\sum_{i=1}^m \bigl|G_{\alpha\mathrm{MSNN}}(x_i^{\!*},y_i^{\!*};\theta)-g_i^{\!*}\bigr|^2.
\end{equation}
This loss penalizes deviations from the coarse-grid reference values during offline training.
Thus, we propose to use the following training loss that combines residual, boundary, and auxiliary terms:
\begin{equation} \label{eq:full_loss}
\mathcal{L}[G_{\alpha\mathrm{MSNN}}] = 
w_{\text{res}}\,\mathcal{L}_{\text{res}} + 
w_{\text{bc}}\,\mathcal{L}_{\text{bc}} + 
w_{\text{aux}}\,\mathcal{L}_{\text{aux}},
\end{equation}
where the boundary enforces the Dirichlet condition:
\begin{equation}
\mathcal{L}_{\text{bc}} \;=\; \frac{1}{|S_{\text{bc}}|} 
\sum_{(x,y)\in S_{\text{bc}}} 
\bigl|G_{\alpha\mathrm{MSNN}}(x,y;\theta)\bigr|^2,
\end{equation}
where $S_{\text{bc}}$ denotes collocation points on the boundary.
The weights $w_{\text{res}}, w_{\text{bc}},$ and $w_{\text{aux}}$ act as balancing hyperparameters and may be adjusted between training stages (see Section~\ref{sec:protocol}).

The auxiliary (anchor) loss provides coarse but reliable information that guides optimization toward physically consistent solutions. Figure~\ref{fig:full_comparison} compares residual-only training with training augmented by coarse-grid anchors for $\mathcal L=-0.1\Delta-50$ on $[0,1]$ at $y=0.65$. The Gaussian source has standard deviation $\sigma=0.01$, corresponding to $\varepsilon=\sqrt{2}\sigma$ in~\eqref{eq:delta}. In the top panel, both models give differential-operator outputs that closely follow this source after the same training budget.
Yet, as the bottom panel reveals, the residual-only network approximates a Green's function slice $G(x,0.65)$ that deviates significantly from the coarse reference.
In contrast, the model trained with anchors accurately reproduces the reference profile and captures the expected oscillatory structure. This comparison illustrates the benefit of optional anchor supervision at the same training budget.

\subsection{Complete Learning Framework and Training Protocol}
\label{sec:protocol}
We now integrate the previously introduced multiscale neural architecture and auxiliary data into a unified learning framework.
The framework is built on three guiding principles.
First, targeted data sampling ensures both physical consistency and identifiability by aligning collocation points and auxiliary reference data with the loss components.
Second, a multi-$\varepsilon$ staged training protocol progressively refines the network across decreasing source widths \rev{$\varepsilon_1 > \varepsilon_2 > \cdots$}, improving stability and enabling accurate learning of small-scale structures.
Third, an overlapping domain decomposition introduces a domain decomposition network $r(y)$ that assigns smooth, nonnegative weights (summing to one) to a small set of shared submodels, allowing local specialization while preserving global continuity in $y$.

\begin{figure}[htbp]
    \centering
    \includegraphics[width=0.95\linewidth]{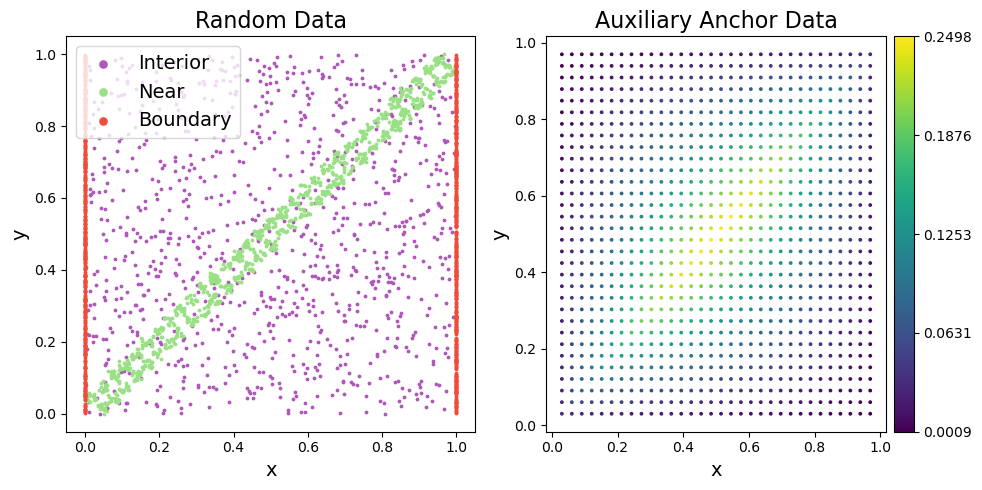}
    \caption{Data sampling strategy for a 1D problem. 
    \textit{Left:} Collocation points include uniformly distributed interior points (purple), denser points near the diagonal $y = x$ (green) to address the singularity, and boundary points (red). 
    \textit{Right:} Auxiliary reference data on a uniform grid, generated from a coarse finite-difference inverse.
    \label{fig:data_sampling}}
\end{figure}

\emph{Training data sampling.}  
The effectiveness of the training depends on a sampling strategy that is carefully aligned with each component of the loss. As shown in Figure~\ref{fig:data_sampling}, we employ two complementary data sources.  
The first source consists of collocation points used to evaluate the PDE residual and enforce boundary conditions. These are drawn from three distributions:  
(a) uniformly distributed interior points to enforce the operator across the domain,  
(b) a denser set near the diagonal $y = x$ to address the singularity of the Green's function, and  
(c) boundary points to impose Dirichlet conditions.  The second source provides auxiliary reference values of the Green's function, obtained by inverting a coarse finite-difference matrix on a uniform grid (Figure~\ref{fig:data_sampling}, right).

\paragraph{A multi-stage training protocol}
{Direct training for a very small $\varepsilon$ is numerically challenging. 
We therefore extend the framework with a coarse-to-fine schedule over the source width.
The NNs are trained in stages of increasing complexity: at each stage, we decrease $\varepsilon$ and augment the active set of branches specialized to the current width, optimizing all active components jointly.
Let $\varepsilon_1>\cdots>\varepsilon_M$ be a prescribed decreasing sequence. 
After $M$ stages, the predictor takes the composite form
\begin{equation}
G_{\alpha\mathrm{MSNN}}(x, y; \theta)  = \sum_{k=1}^M\sum_{j=1}^{L-1} 
\phi_{j,k}^{\varepsilon_k}(x,y;\alpha_{j,k}(y),\beta_{j,k}(y),\theta_{j,k}) + \phi_{\text{far}}(x,y;\theta_{\text{far}}),
\end{equation}
where $j$ indexes the $L-1$ near- and middle-field branches and $k$ indexes the source width.
Each local branch takes the form
$\phi_{j,k}^{\varepsilon_k}=\varepsilon_k^{\alpha_{j,k}(y)} 
\hat{f}_{j,k}\!\big((x-y)/\varepsilon_k^{\beta_{j,k}(y)},y;\,\theta_{j,k}\big)$
with a smooth saturating activation in $f_{j,k}$; the global far-field is shared across stages.
Within each stage we impose the ordering $\beta_{j+1,k}<\beta_{j,k}$ and $\alpha_{j+1,k}\!\ge\!\alpha_{j,k}$ so that support broadens while amplitude weakens across $j$.
During stage $k$, only the partial sum with $k'\le k$ and the far-field $\phi_{\text{far}}$ shared across stages is active; when moving to $\varepsilon_{k+1}$ we add the new near-field branch with $\beta_{1,k+1}=1$, then fine-tune jointly.
At the transition $k\!\to\!k{+}1$ the right-hand side also change from $\mathcal{N}_{\varepsilon_k}$ to $\mathcal{N}_{\varepsilon_{k+1}}$.

We further enhanced the stability of the training by dynamically adjusting key hyperparameters. The auxiliary loss weight $w_{\text{aux}}$ is set large at the start (e.g., $w_{\text{aux}}=1.0$) to guide the network toward the correct global structure using coarse grid anchor points. As training proceeds, $w_{\text{aux}}$ is gradually decreased so that optimization focuses more on minimizing the PDE residual and enforcing physical constraints.

\iffalse
Regarding optimizers, first-order methods such as Adam are well-suited to rapid exploration: they quickly reduce the residual loss and guide the model toward a promising basin of attraction. However, Adam alone often stalls before the network converges to the precise $(\alpha_j,\beta_j,\theta_j)$ configuration needed for accurate multiscale balance. In contrast, second-order methods can refine parameters once a good basin is reached, but are too expensive in the early stages. 
We use Adam in all previous phases, and transition to second-order methods during the final joint training of the entire model when Adam converges slowly.
\fi 

\paragraph{An overlapping domain decomposition (DD) training approach} Global NNs can be insufficient for challenging problems with small $\varepsilon$, especially when the Green's function varies significantly.
Since the boundary condition for Green's function only enforces on $x$, it is natural to apply DD, and learn multiple approximations in the form $G_\varepsilon(x,y),x\in\Omega,y\in\Omega_i$ with $\bigcup_i\Omega_i=\Omega$ \cite{teng2022learning,hao2024multiscale}.
However, a naive partition can impose boundaries misaligned with the underlying physics and, more critically, create non-physical discontinuities in the learned Green's function.

Motivated by the  Mixture-of-Experts (MoE) approach  \cite{shazeer2017outrageously}, we introduce a source-conditioned overlapping 
DD approach,
avoiding the drawbacks of static geometric partitioning. 
This approach is activated only as a final specialization stage after all $M$ levels are trained. We begin from the fully trained global model of the last level $\sum_{j=1}^{L-1}\phi_{j,M}$. We then replicate this model into $Q$ models and add small parameter perturbations to break symmetry. In parallel, we initialize a lightweight DD network $r(y)$ with a softmax output so that $r_m(y)\ge 0$ and $\sum_m r_m(y)=1$, and pretrain it to assign near-unity weight on a few seed points for each subdomain, as illustrated for a 2D problem in Figure~\ref{fig:domain decomposition}. The predictor is then switched to the soft mixture
\begin{align}
G_{\alpha\mathrm{MSNN}}(x, y; \theta)  &= \sum_{m=1}^{Q} r_m(y)\sum_{j=1}^{L-1} 
\phi_{j,M,m}(x,y;\alpha_{j,M,m}(y),\beta_{j,M,m}(y),\theta_{j,M,m},\varepsilon_M) \\ &+ \sum_{k=1}^{M-1}\sum_{j=1}^{L-1} 
\phi_{j,k}(x,y;\alpha_{j,k}(y),\beta_{j,k}(y),\theta_{j,k},\varepsilon_k) + \phi_{\text{far}}(x,y;\theta_{\text{far}}),
\end{align}

and training continues jointly, letting the optimization discover an effective overlapping DD without graph partitioning. 
Note that instead of the top-$k$ strategy used in the standard MoE approach, we do a full weighted sum to ensure the smoothness of the learned network.
Because the regions where $r_m(y)$ are positive overlap, the composite $G(x_0,\cdot)$ remains smooth in $y$. 
This smoothness is important for the integration over $y$ and, as discussed next, for the compressibility of the learned dense operator.

\begin{figure}[htbp]
    \centering
    \includegraphics[width=0.95\linewidth]{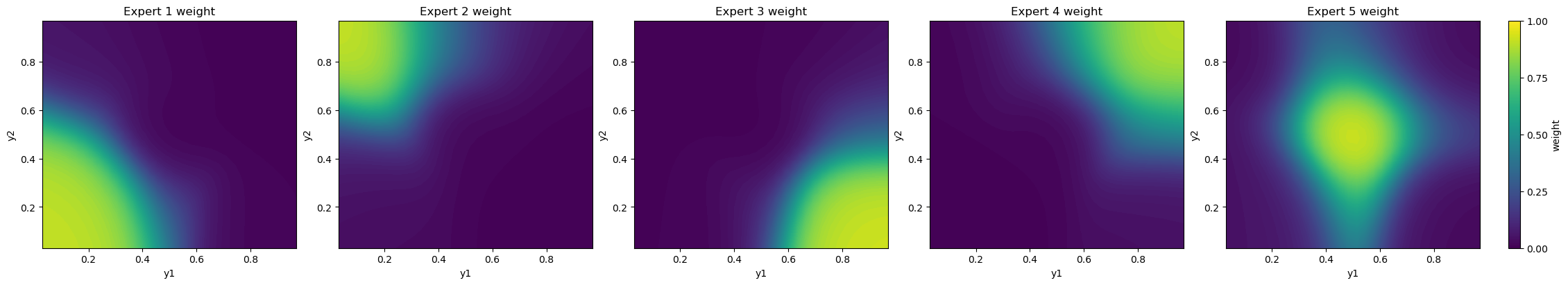}
    \caption{Initial weights of the overlapping DD network on the source domain $y$ in $[0,1]^2$. The network is pre-trained to softly target the center and four corner regions, providing a structured initialization before joint fine-tuning.}
    \label{fig:domain decomposition}
\end{figure}

\begin{figure}[htbp]
    \centering
    \includegraphics[width=0.65\linewidth]{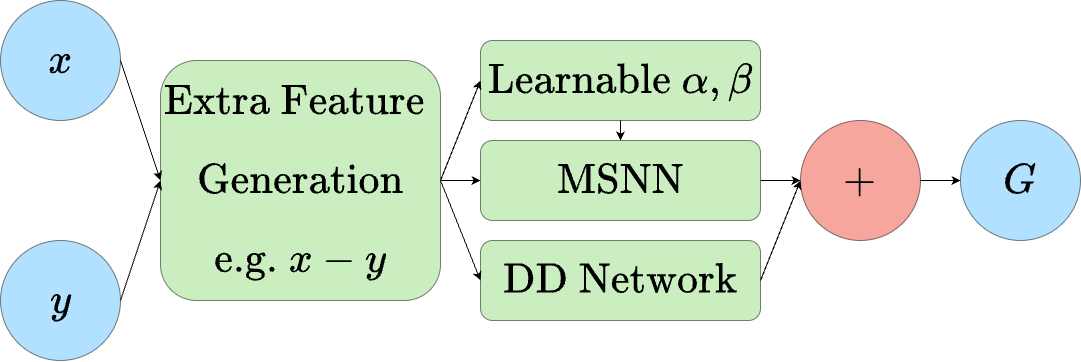}
    \caption{The design of the complete $\alpha\text{MSNN}$ framework.}
    \label{fig:diagram}
\end{figure}

Finally, in Figure~\ref{fig:diagram} 
we summarize
the complete $\alpha$MSNN pipeline, showing how the learnable scaling exponents $(\alpha,\beta)$, the multiscale backbone, and the overlapping DD combine to produce the Green's function $G(x,y)$ from inputs $x$ and $y$.}
\rev{Algorithm 3.1 summarizes the training procedure, where \(T_k\) and \(T_{\mathrm{MoE}}\) denote the numbers of optimization steps in the corresponding stages.}

\begin{algorithm}[htbp]

\caption{\rev{Multistage training of the $\alpha$MSNN Green kernel.}}
\label{alg:green_training}
\begin{algorithmic}[1]
\STATE Initialize the shared far-field network; prepare coarse anchor data if enabled.
\FOR{$k=1,\ldots,M$}
  \STATE Add branches at width $\varepsilon_k$, retaining earlier branches and their source widths.
  \FOR{$t=1,\ldots,T_k$}
    \STATE Sample collocation points and anchor data.
    \STATE Sample operator parameters for parametric problems.
    \STATE Jointly update the active parameters using loss~\eqref{eq:full_loss}.
  \ENDFOR
\ENDFOR
\IF{MoE refinement is enabled}
  \STATE Replace the finest-scale component with $Q$ perturbed expert copies.
  \STATE Initialize and pretrain the gate $r(y)$.
  \STATE Jointly train all branches and $r(y)$ for $T_{\mathrm{MoE}}$ steps at $\varepsilon_M$.
\ENDIF
\STATE \textbf{return} the trained Green kernel $G_{\alpha\mathrm{MSNN}}$.
\end{algorithmic}
\end{algorithm}

\rev{In practice, training can start with a relatively large source
width $\varepsilon$ that the collocation samples can resolve.
The initial network size can be guided by experience with PINNs
for similar PDEs, since our training also minimizes PDE and
boundary residuals. New branches are then added as $\varepsilon$
decreases, with the $\varepsilon$-dependence of the network-size
estimates in~\cite{hao2024multiscale} providing guidance for
increasing the network capacity. PDE and boundary residuals on
independent validation samples, together with a few trial solves,
can guide the choice of network size and the number of branches.}

% Original submission section: matrix.tex
\section{Converting Green's function into a Data-sparse format}\label{sec:preconditioner}
Once the Green's function has been learned, our goal is to construct a practical and efficient preconditioner
$\mathbf{M}^{-1} \;\approx\; \mathbf{G}$,
that balances three key objectives: (i) low memory usage, (ii) fast matrix--vector multiplication, and (iii) sufficient spectral quality to accelerate iterative solvers. While the dense matrix $\mathbf{G}$ formed from $G_{\alpha\mathrm{MSNN}}$ provides an accurate surrogate for $\mathbf{A}^{-1}$, directly storing or applying this operator is infeasible at large scale.  

\subsection[Sparse and hierarchical representations]{\rev{Sparse and hierarchical representations}}
\label{sec:data-sparse}

A dense discretization of the learned Green's function provides a direct approximation to the inverse matrix, but its storage and application costs scale quadratically with the number of degrees of freedom. To obtain more efficient representations, we exploit the spatial structure of the learned kernel. We consider two complementary cases: hierarchical compression when significant long-range interactions are present, and sparse representations when the Green's function is strongly localized.

\rev{For kernels with significant long-range interactions, spatial clustering naturally leads to a hierarchical matrix (\(\mathcal H\)-matrix) representation, in which near-field interactions are stored in dense blocks while admissible far-field blocks are approximated by low-rank matrices~\cite{hbook,ddh2}. Under the usual assumptions of uniformly bounded block ranks and a suitably balanced block tree, the resulting representation can achieve near-linear storage and application complexity. These estimates concern only the compressed operator and do not include the cost of learning the Green's function or constructing the hierarchical representation.
When the learned Green's function is strongly localized, a simpler alternative is to retain only selected entries of the learned inverse matrix \(\mathbf G\approx\mathbf A^{-1}\), yielding a sparse approximate inverse that can be applied directly by sparse matrix-vector multiplication.}

\subsection[Schur complement preconditioning]{\rev{Schur complement preconditioning}}
\label{sec:schur_deployment}
{

The learned Green's function also fits naturally into a two-level ILU framework based on domain decomposition~\cite{li_saad_sosonkina_2003_parms,gpuilu}.
We can partition the physical domain into $p$ nonoverlapping subdomains separated by an interface $\Gamma$, as illustrated in the left panel of Figure~\ref{fig:green_schur_dd}.

\begin{figure}[htbp]
  
  \centering
  \includegraphics[width=\textwidth]{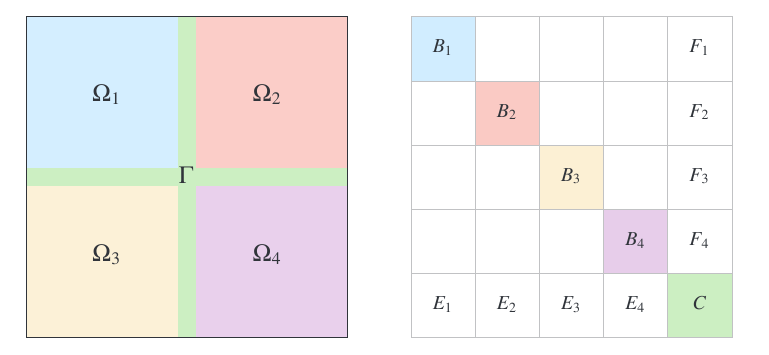}
  \caption{\rev{Domain decomposition and the corresponding matrix structure.
  \textit{Left:} A nonoverlapping decomposition of the physical domain with separator $\Gamma$.
  \textit{Right:} The matrix after ordering the subdomain interiors before the separator.}}
  \label{fig:green_schur_dd}
\end{figure}

Ordering the subdomain interiors before the separator with permutation matrix $\mathbf P$ gives the block structure in Figure~\ref{fig:green_schur_dd} (right).
For nonsingular interior blocks, the reordered matrix and its Schur complement $\mathbf S$ are
\begin{equation}\label{eq:schur_blocks}
  \mathbf A_0=\mathbf P^T\mathbf A\mathbf P
  =\begin{bmatrix}\mathbf B&\mathbf F\\\mathbf E&\mathbf C\end{bmatrix},
  \qquad \mathbf S=\mathbf C-\mathbf E\mathbf B^{-1}\mathbf F.
\end{equation}
Here $\mathbf B=\operatorname{diag}(\mathbf B_1,\ldots,\mathbf B_p)$ contains the independent interior blocks.

For nonsingular $\mathbf A$, the inverse-block identity gives a direct connection between the learned Green's function and the Schur complement:
\begin{equation}\label{eq:green_schur_inverse}
  (\mathbf A_0^{-1})_{\Gamma,\Gamma}=\mathbf S^{-1},
  \qquad
  \mathbf T=(\mathbf P^T\mathbf G\mathbf P)_{\Gamma,\Gamma}
  \approx\mathbf S^{-1},
\end{equation}
where $\mathbf G$ denotes the normalized learned discrete inverse in the original ordering. Thus, only kernel values associated with separator coordinate pairs are needed to construct $\mathbf T$; the full learned inverse matrix $\mathbf G$ need not be assembled. If the separator contains $n_\Gamma$ unknowns, constructing $\mathbf T$ requires $n_\Gamma^2$ kernel evaluations, while storing the resulting dense matrix and applying it to a vector each require $\mathcal O(n_\Gamma^2)$ work or storage. The interior blocks in $\mathbf{B}$ are treated separately using local incomplete factorizations.

For this construction, training follows Algorithm~\ref{alg:green_training}, with source points concentrated near the separator $\Gamma$, target points sampled throughout $\Omega$, and boundary points on $\partial\Omega$. The resulting approximation $\mathbf T\approx\mathbf S^{-1}$ is then combined with local ILUT factorizations
\[
\mathbf L_B\mathbf U_B\approx\mathbf B.
\]
Defining
\[
\mathbf M_B=\mathbf U_B^{-1}\mathbf L_B^{-1}\approx\mathbf B^{-1},
\]
we obtain a two-level ILU preconditioner in which $\mathbf M_B$ provides the local interior solves and $\mathbf T$ supplies an approximate inverse on the separator. Algorithm~\ref{alg:green_schur} summarizes the resulting setup and application.

\begin{algorithm}[htbp]

\caption{\rev{Two-level ILU preconditioning with a learned Green's function.}}
\label{alg:green_schur}
\begin{algorithmic}[1]
\STATE \textbf{Setup:} given $\mathbf A_0$, $\mathbf P$, and a trained Green function $G$.
\STATE Compute local ILUT factors $\mathbf L_B\mathbf U_B\approx\mathbf B$ and set $\mathbf M_B=\mathbf U_B^{-1}\mathbf L_B^{-1}$.
\STATE Evaluate and normalize separator kernel values to assemble $\mathbf T\approx\mathbf S^{-1}$.
\STATE \textbf{Apply:} given the residual $\mathbf r$ in the original ordering.
\STATE Set $[\mathbf r_B^T,\mathbf r_\Gamma^T]^T=\mathbf P^T\mathbf r$.
\STATE Compute $\mathbf u_B=\mathbf M_B\mathbf r_B$ and $\mathbf g_\Gamma=\mathbf r_\Gamma-\mathbf E\mathbf u_B$.
\STATE Compute $\mathbf z_\Gamma=\mathbf T\mathbf g_\Gamma$.
\STATE Recover $\mathbf z_B=\mathbf M_B(\mathbf r_B-\mathbf F\mathbf z_\Gamma)$.
\STATE Return $\mathbf z=\mathbf P[\mathbf z_B^T,\mathbf z_\Gamma^T]^T$.
\end{algorithmic}
\end{algorithm}

Constructing $\mathbf T$ directly from the learned Green's function avoids forming the Schur complement through interior elimination and does not depend on the local ILUT factors. On a distributed-memory system, each process can independently evaluate and scale its assigned entries of $\mathbf T$ without communication, once the trained network and the required interface coordinates and weights are locally available. The local ILUT factorizations can be performed independently across subdomains using existing CPU implementations, while
$\mathbf T$ is constructed concurrently on GPUs.
}

% Original submission section: exp.tex

\section{Numerical results}\label{sec:experiments}
This section evaluates the proposed NN approach for approximating Green’s functions 
and constructing the corresponding inverse preconditioner, represented either as a dense matrix 
or in compressed $\mathcal{H}$-matrix form. 
\rev{The first three experiments examine the convergence behavior and spectral properties of the learned Green kernel used as a standalone preconditioner. The final experiment evaluates its practical performance within a two-level ILU framework, comparing setup and solve times with classical two-level ILUT and two-level Schwarz--ILUT.}

Because the coefficient matrices are generally nonsymmetric and indefinite, 
we employ the Flexible GMRES (FGMRES) algorithm as the Krylov accelerator. 
In the first three experiments, we compare three configurations across a suite of linear systems:
(i) unpreconditioned FGMRES,
(ii) FGMRES preconditioned with the dense $\alpha$MSNN Green’s kernel, and
(iii) FGMRES preconditioned with its $\mathcal{H}$-matrix compression. 
Convergence is measured by the number of iterations required for the relative residual 
to fall below $10^{-6}$; runs that do not reach this tolerance within the iteration cap of $500$ 
are denoted by “F.” 
For each configuration, three reference right-hand sides with i.i.d.\ entries in $[-0.5,0.5]$ 
are used, starting from zero initial vectors. 
We report both the mean and standard deviation of iteration counts over these runs. 
All experiments use FGMRES with a restart value of~50 (i.e., FGMRES(50)).

\rev{The experiments were conducted on a laptop equipped with
an Intel Core Ultra 9 275HX processor (24 CPU cores: 8 performance cores and
16 efficiency cores), an NVIDIA GeForce RTX 5090 Laptop GPU with 24~GB of
dedicated GPU memory, and 32~GB of system memory.
The software environment comprised Ubuntu~24.04.4 LTS running under WSL2 on
Windows~11, Python~3.12.3, JAX~(v0.11.1), Flax~(v0.12.9), Optax~(v0.2.8),
NumPy~(v2.4.4), and SciPy~(v1.17.1).
WSL2 exposed 20 logical processors and approximately 23.5~GiB of memory.
Training and Green-kernel evaluation used GPU single-precision arithmetic,
while linear solves and external residual checks used double precision.
BLAS and OpenMP threads were limited to one.
The GPU training times reported below cover the main Adam phases, including JIT compilation.}

\subsection{Model Problem}

Throughout the entire experiments section, we consider the following model problem with a Dirichlet boundary condition:
\begin{align}
-\nabla\cdot\big(a(x)\nabla u(x)\big) + \mathbf{b}(x)\cdot\nabla u(x) + c(x)u(x) &= f(x)\ \ \text{in} \ \Omega, \nonumber \\
u(x) &= 0 \ \ \text{on} \ \partial \Omega,
\end{align}
where the problem is defined on the unit cube domain $\Omega=[0,1]^d$ for $d=1, 2$.
Discretization is performed on a uniform grid using a central difference scheme, except in convection-dominated cases where an upwind scheme is applied to the convective term for stability.
We test convection-dominated and reaction--diffusion problems with diagonal $a(x)$ (five-point stencils in 2D and three-point stencils in 1D), and anisotropic rotated diffusion problems with a nine-point stencil. The reaction--diffusion tests include indefinite operators.

\subsection{Experiment 1: 1D Problems}

We begin with one-dimensional test problems and construct two cases by varying the PDE coefficients: 
a convection-dominated case with $a(x)=0.01$, $b(x)=1+x^2$, $c(x)=0$, 
and an indefinite reaction--diffusion case with $a(x)=1$, $b(x)=0$, $c(x)=-50(1+x^2)$. 
These configurations yield linear systems with distinct characteristics: 
the convection-dominated case produces a highly nonsymmetric matrix, 
while the second case has both positive and negative eigenvalues.

\paragraph{Network architecture}
We employ a 3-scale $\alpha$MSNN with $\varepsilon\in\{10^{-2},\,10^{-3},\,10^{-4}\}$. 
The architecture comprises one global far-field branch and, for each $\varepsilon$, 
a near-field and a middle-field branch. 
{All networks have three hidden layers.}
The far-field network has $50$ neurons per layer. 
For the near- and middle-field branches at $\varepsilon=10^{-2}$, $10^{-3}$, and $10^{-4}$, 
we use $10$, $15$, and $20$ neurons per hidden layer, respectively. 
At the smallest scale ($\varepsilon=10^{-4}$), 
the overlapping domain decomposition (DD) strategy is activated. 
We partition the domain into three subdomains centered at $y\in\{0.0,\,0.5,\,1.0\}$ 
and deploy a lightweight DD network together with scaling networks for 
$\alpha(y)$ and $\beta(y)$. 
Each of these auxiliary networks uses $5$ neurons per hidden layer.  %
All networks employ $\tanh$ activations. 
For deep architectures, the $\beta$ scaling factor is applied only to the input layer, 
while the $\alpha$ scaling factor is applied only to the output layer.

\paragraph{Training configuration}
Training is performed with the Adam optimizer. 
We use $500$ epochs for the two larger $\varepsilon$ levels, 
$1000$ epochs to initialize the smallest scale ($\varepsilon=10^{-4}$), 
and an additional $3000$ epochs for the overlapping DD specialization. 
Each epoch consists of $20$ mini-batches. 
\rev{Training took 18.05 and 17.91 minutes for the convection-dominated and indefinite cases, respectively.}
Each mini-batch contains:
\begin{itemize}[leftmargin=2em]
    \item $500$ boundary samples enforcing Dirichlet conditions,
    \item {$500$ auxiliary  subsamples from $1024$ anchor data generated using a coarse grid with $34$ points,} %
    \item $500$ uniformly distributed residual pairs over $\Omega\times\Omega$, and
    \item $1500$ near-diagonal residual pairs with $|x-y|\le3\varepsilon$, where $\varepsilon$ corresponds to the current training stage.
\end{itemize}
Once the $\alpha$MSNN is trained, the learned continuous kernel can be discretized on demand 
to construct a preconditioner for any target grid resolution. 
To assess scalability, we evaluate the resulting preconditioners on a sequence of uniform grids 
with $N+2$ points, where $N=8,16,32,\ldots,2048$.
The compressed $\mathcal{H}$-matrix preconditioner is generated directly from the learned 
$\alpha$MSNN kernel. 
The $\mathcal{H}$-matrix compression follows a standard adaptive block-clustering procedure 
with a maximum leaf size of $\min\{\sqrt{N},\,128\}$ 
and an admissibility parameter $\eta=1.0$ to define the block partitioning. 
For each admissible block identified by this hierarchy, 
we first compute an intermediate rank-3 approximation via the Nyström method 
using the nearest interpolation points,
which is then truncated to a final rank-1 representation, 
achieving a compression ratio of up to $46.5\times$ for the largest problem size.

\begin{figure}[htbp]
    \centering
    \includegraphics[width=0.95\linewidth]{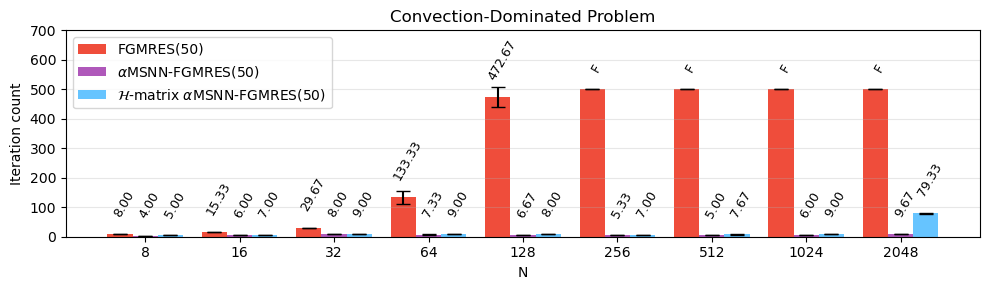}
    \includegraphics[width=0.95\linewidth]{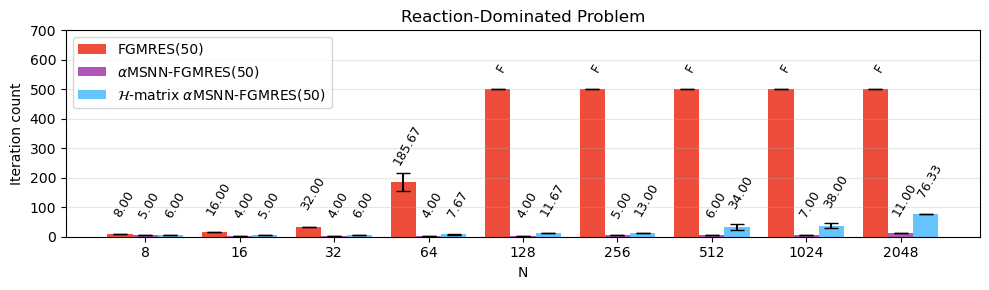}
    \caption{Problem size vs number of FGMRES(50) iterations for two problems: a convection-dominated problem with $a(x)=0.01$, $b(x)=1+x^2$, and $c(x)=0.0$; and an indefinite problem with $a(x)=1.0$, $b(x)=0.0$, and $c(x)=-50(1+x^2)$. F indicates that the solver fails to converge within $500$ iterations.}
    \label{fig:iteration oned}
\end{figure}

Figure~\ref{fig:iteration oned} summarizes the average FGMRES(50) iteration
number
as a function of problem size for the convection-dominated and indefinite cases. 
Error bars indicate the standard deviation over three runs with different random right-hand sides. 
For the unpreconditioned solver, the iteration count increases rapidly with $N$, 
and convergence eventually fails within the 500-iteration limit. 
In contrast, preconditioning with the dense $\alpha$MSNN Green’s kernel 
exhibits excellent scalability: the iteration count remains low and nearly constant 
across all grid sizes for both test problems. 
\rev{Rank-one $\mathcal H$-matrix compression preserves convergence but can substantially increase iteration counts on larger grids ($N=2048$ in the convection-dominated case and $N\geq512$ in the indefinite case). All far-field blocks use the same truncation rank~1 in this experiment.}

\begin{figure}[htbp]
    \centering
    \includegraphics[width=0.95\linewidth]{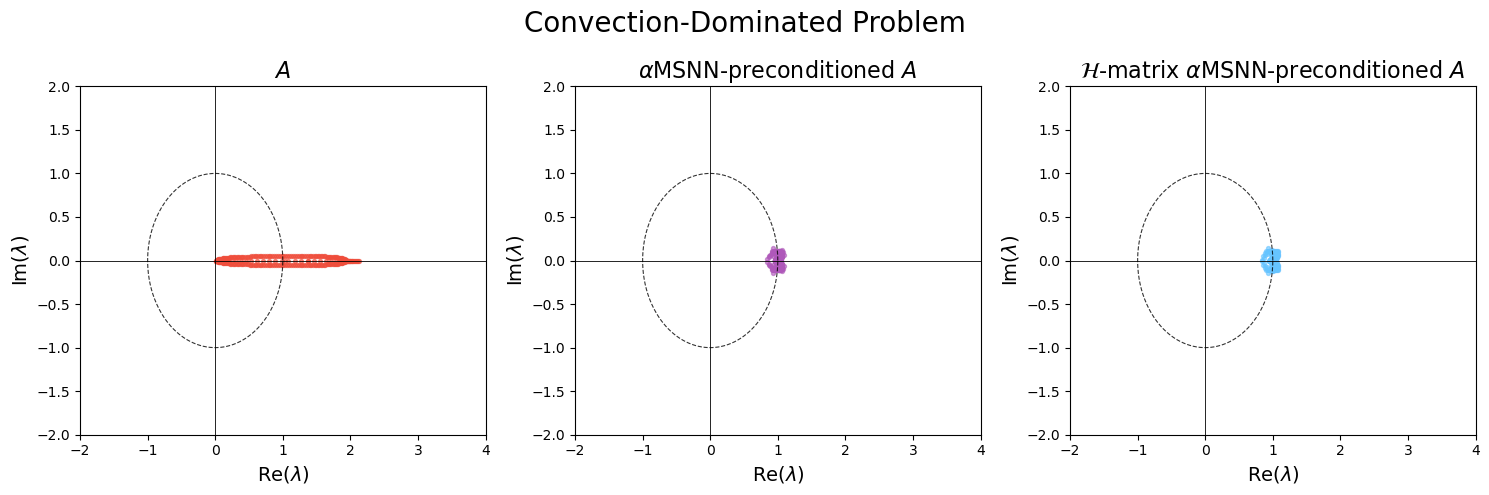}
    \includegraphics[width=0.95\linewidth]{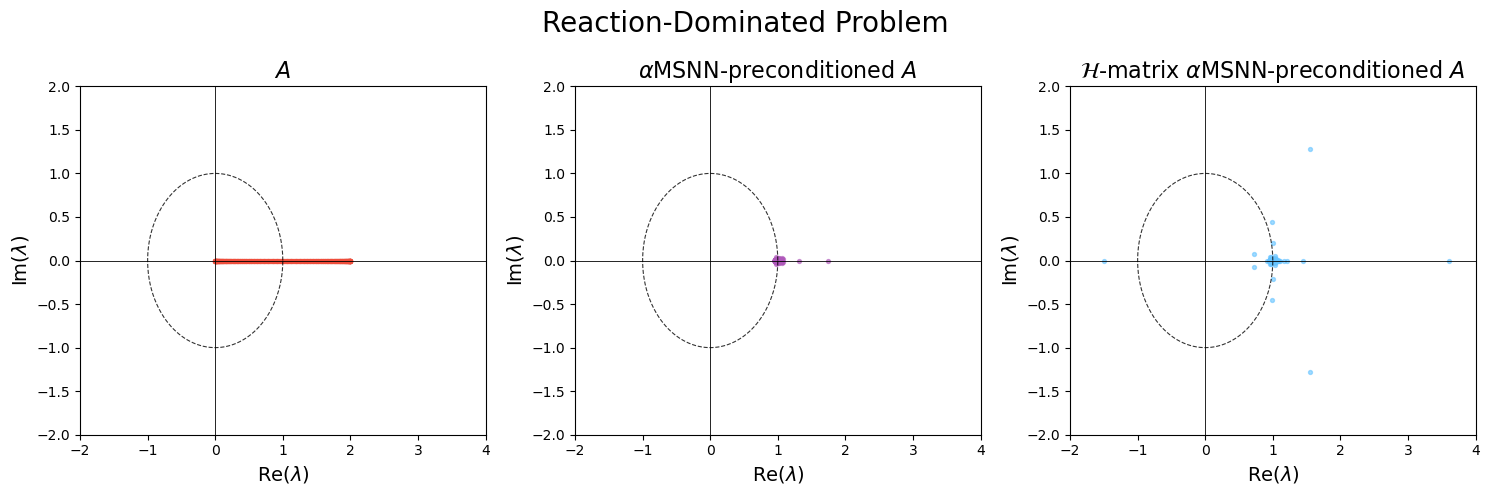}
    \caption{Spectra of the original matrix (left), the preconditioned matrix using the dense $\alpha$MSNN kernel (middle), 
and the preconditioned matrix using the $\mathcal{H}$-matrix–compressed $\alpha$MSNN kernel (right) 
for two problems: a convection-dominated case with $a(x)=0.01$, $b(x)=1+x^2$, $c(x)=0$, 
and an indefinite case with $a(x)=1.0$, $b(x)=0$, $c(x)=-50(1+x^2)$. 
The matrix size is $256\times256$. 
For visual comparison, each spectrum is linearly scaled by the inverse of its mean eigenvalue.
}
    \label{fig:spectrum oned}
\end{figure}

To provide a spectral interpretation of the convergence behavior, 
Figure~\ref{fig:spectrum oned} visualizes the eigenvalue distributions of the original and preconditioned matrices 
for the $N=256$ case. 
The spectra of the original unpreconditioned matrices (left column) highlight their challenging nature: 
the convection-dominated operator exhibits a non-normal, widely spread spectrum, 
while the indefinite operator has eigenvalues distributed on both sides of the imaginary axis. 
These unfavorable spectral properties explain the poor performance of the unpreconditioned solver. 

In contrast, preconditioning with the dense $\alpha$MSNN Green’s kernel (middle column) 
dramatically improves the spectrum, clustering the eigenvalues of both problems 
into a compact region away from the origin.  
The $\mathcal{H}$-matrix version (right column) yields nearly identical results: 
although rank-1 truncation introduces minor perturbations that slightly broaden the cluster, 
all eigenvalues remain well separated from the origin, 
with most confined within a small disk. 
These observations confirm that $\mathcal{H}$-matrix compression preserves 
the essential spectral properties of the dense preconditioner, 
accounting for their comparable convergence behavior.

Finally, note that symmetry was not explicitly enforced on the preconditioner in these experiments, 
since FGMRES does not require it. 
If a symmetric preconditioner is desired, it can be readily imposed 
by augmenting the training loss with a symmetry penalty 
$|\phi(x,y) - \phi(y,x)|$~\cite{teng2022learning} 
and constructing the final operator from the symmetrized kernel 
$(\phi(x,y) + \phi(y,x))/2$.

\subsection{Experiment 2: 2D Problems}

We next extend our evaluation to two-dimensional problems on the unit square $[0,1]^2$. 
Following the same methodology, we construct two test cases that are particularly challenging for computing effective preconditioners. 
The first is a convection-dominated problem with 
\begin{align}
a(x_1,x_2)&=\operatorname{diag}\{0.01(1+x_1^2+x_2^2),\,0.01(1+x_1^2+x_2^2)\},\\
b(x_1,x_2)&=[\,1+x_2^2,\,1+x_1^2\,],\quad
c(x_1,x_2)=0.
\end{align}
The second is a reaction-dominated problem with 
\begin{equation}
a(x_1,x_2)=\operatorname{diag}\{1,1\},\quad
b(x_1,x_2)=[\,0,0\,],\quad
c(x_1,x_2)=-10(1+x_1^2+x_2^2).
\end{equation}

The $\alpha$MSNN architecture, including the number of hidden layers and neurons per layer, remains identical to the 1D setup. However, to account for the increased problem complexity, we adopt a more conservative $\varepsilon$-schedule with $\varepsilon\in\{10^{-2},5\times10^{-3},2\times10^{-3}\}$.
The framework employs five subdomains with one centered at $(0.5,0.5)$ and four targeted to the domain corners.
The training epochs are adjusted to $1000$ epochs for the top levels, $1000$ epochs for initializing the lowest level, and $2000$ epochs for domain decomposition.
\rev{A four-stage GPU run took 21.33 minutes for the convection-dominated model.}
All other experimental parameters, including the data sampling strategy, are kept the same as in the 1D case, with the key exception that the auxiliary anchor data is now generated from a 2D $34\times34$ grid.

For the 2D experiments, we discretize the trained $\alpha$MSNN on a sequence of uniform $(N+2)\times(N+2)$ grids for $N=32,48,64,96,128$.
To handle the increased complexity of the 2D problems, the $\mathcal{H}$-matrix construction employs a more sophisticated strategy.
The block clustering is configured with an admissibility parameter $\eta=0.7$ and a maximum leaf size $n_l$ of $64$ for $N\leq64$ and $128$ for $N>64$.
The target rank $k$ for each admissible block is determined by a size-dependent rank rule based on both the global problem size $N$ and the local block size $n_b$, defined as the smallest between number of rows and number of columns. 
First, a base rank is established according to the overall problem size: $k=10$ for problems with $N\leq64$, which is increased to $k=20$ for the larger problems where $N>64$. When $n_b$ is larger than maximum leaf size $n_l$, this base rank is then scaled up logarithmically for any admissible block of size $n$ by a factor of $\log(n_b)/\log(n_l)$. All low-rank approximations are computed using an oversampled random Nystr\"om method, where an initial approximation at twice the target rank is computed and then truncated.

\begin{figure}[htbp]
    \centering
    \includegraphics[width=0.95\linewidth]{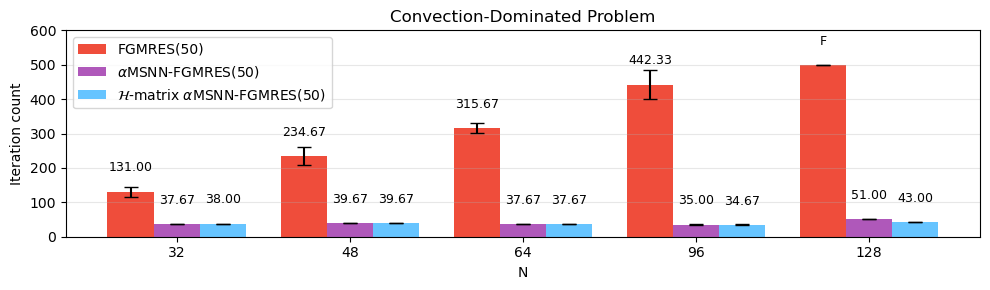}
    \includegraphics[width=0.95\linewidth]{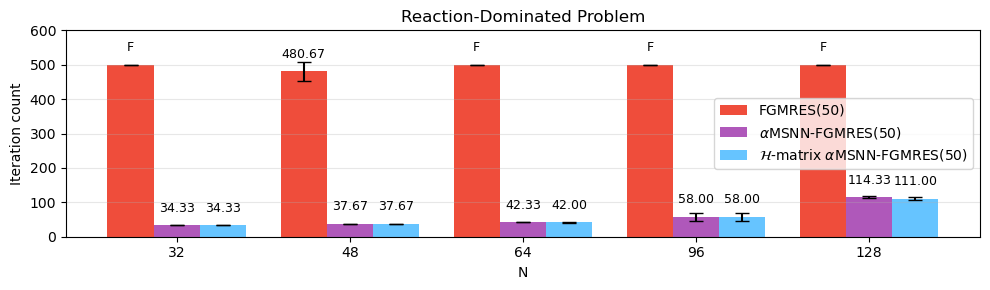}
    \caption{Problem size vs number of FGMRES(50) iterations for two problems: a convection-dominated problem with $a(x_1,x_2)=\text{diag}\{0.01[1+x_1^2+x_2^2,1+x_1^2+x_2^2]\}$, $b(x_1,x_2)=[1+x_2^2,1+x_1^2]$, and $c(x_1,x_2)=0.0$; and a reaction-dominated problem with $a(x_1,x_2)=\text{diag}\{[1.0,1.0]\}$, $b(x_1,x_2)=[0.0,0.0]$, and $c(x_1,x_2)=-10(1+x_1^2+x_2^2)$. F indicates that the solver fails to converge within $500$ iterations.}
    \label{fig:iteration twod}
\end{figure}

Figure~\ref{fig:iteration twod} summarizes the results for the 2-D problems. 
As in the 1-D case, the unpreconditioned solver deteriorates rapidly with increasing $N$ 
and often fails to reach the prescribed tolerance within the iteration budget. 
In contrast, the dense $\alpha$MSNN preconditioner exhibits good scalability: the iteration count increases only mildly with grid size for both the convection- and reaction-dominated cases, with a slightly larger rise at $N=128$ in the reaction-dominated problem.
The $\mathcal{H}$-matrix variant closely follows the dense baseline, 
and in some instances even shows marginally better convergence.
It is worth noting that for the $128\times128$ case, the $\mathcal{H}$-matrix preconditioner achieves 
approximately a $4\times$ reduction in storage compared with the dense $\alpha$MSNN, 
while maintaining nearly identical iteration counts. 
This result was obtained using a deliberately conservative setup, 
employing fixed-rank Nystr\"om approximation, 
and no nested bases. 
A higher compression ratio and faster application 
are expected with $\mathcal{H}^2$ formats, which can further enhance 
the overall performance.

\subsection{Experiment 3: Parameterized Preconditioner}
\label{subsec:2d_anisotropic}
In this experiment, we extend the approach from a single PDE to a family of parameterized problems, training a single model to serve as a preconditioner across the entire set. 
The key idea is that $\alpha$MSNN decomposes Green’s function into near-, middle-, and far-field components with learnable scales and DD, making each component easier to approximate and stabilizing training for a fixed operator. 
Building on this foundation, we generalize $\alpha$MSNN to learn a unified model for parameterized PDEs, with the goal of producing effective preconditioners that adapt to varying coefficients within the same problem family.
For this study, we focus on a pure diffusion problem corresponding to the 2-D anisotropic rotated Laplacian with Dirichlet boundary conditions, where both convection and reaction terms are set to zero. 
The anisotropic diffusion tensor is parameterized by an anisotropy ratio $\xi$ and rotation angle $\theta$, defined as:

\begin{equation}
a(x) =
    \begin{bmatrix}
        \cos^2 \theta + \xi \sin^2 \theta & \cos \theta \sin \theta (1 - \xi) \\
        \cos \theta \sin \theta (1 - \xi) & \sin^2 \theta + \xi \cos^2 \theta
    \end{bmatrix}.
    \label{eq:rlap}
\end{equation}
We fix $\xi = 0.1$ and aim to learn a parametric Green’s function 
$G(x, y, \theta)$ corresponding to the PDE operator. 
Since the coefficient is $\pi$-periodic in $\theta$, 
we encode the parameter using $(\cos 2\theta,\, \sin 2\theta)$ 
and train $\alpha$MSNN with input $(x, y, \cos 2\theta, \sin 2\theta)$. 
Remarkably, this extension requires no modification to the overall framework 
other than providing the two additional input features. 
The $\beta$ scaling acts only on the spatial coordinates $(x, y)$, 
while the parameter channels $(\cos 2\theta,\, \sin 2\theta)$ are passed through unscaled.

We employ the same three-scale $\alpha$MSNN architecture as in the previous 2D experiments, 
with all hidden-layer widths doubled to account for the increased complexity. 
The DD framework is configured with nine overlapping subdomains 
in the angular parameter $\theta$, effectively partitioning the parameter space. 
During training, each mini-batch samples not only the spatial coordinates 
but also the rotation angle $\theta$. 
To accommodate this additional dimension, 
we double the number of points per mini-batch 
while keeping the proportions of boundary, anchor, and interior samples unchanged. 
The auxiliary anchor data are similarly extended: 
they are generated on the same $34\times 34$ spatial grid 
but include four representative angles 
$\theta \in \{0,\, \pi/4,\, \pi/2,\, 3\pi/4\}$. 
Training uses the Adam optimizer for 1000 epochs per stage, totaling 4000 epochs. 
\rev{Training took 23.93 minutes on the GPU.}
The $\mathcal{H}$-matrix construction parameters and adaptive strategy 
are identical to those described in the previous section.

\begin{figure}
    \centering
    \includegraphics[width=0.95\linewidth]{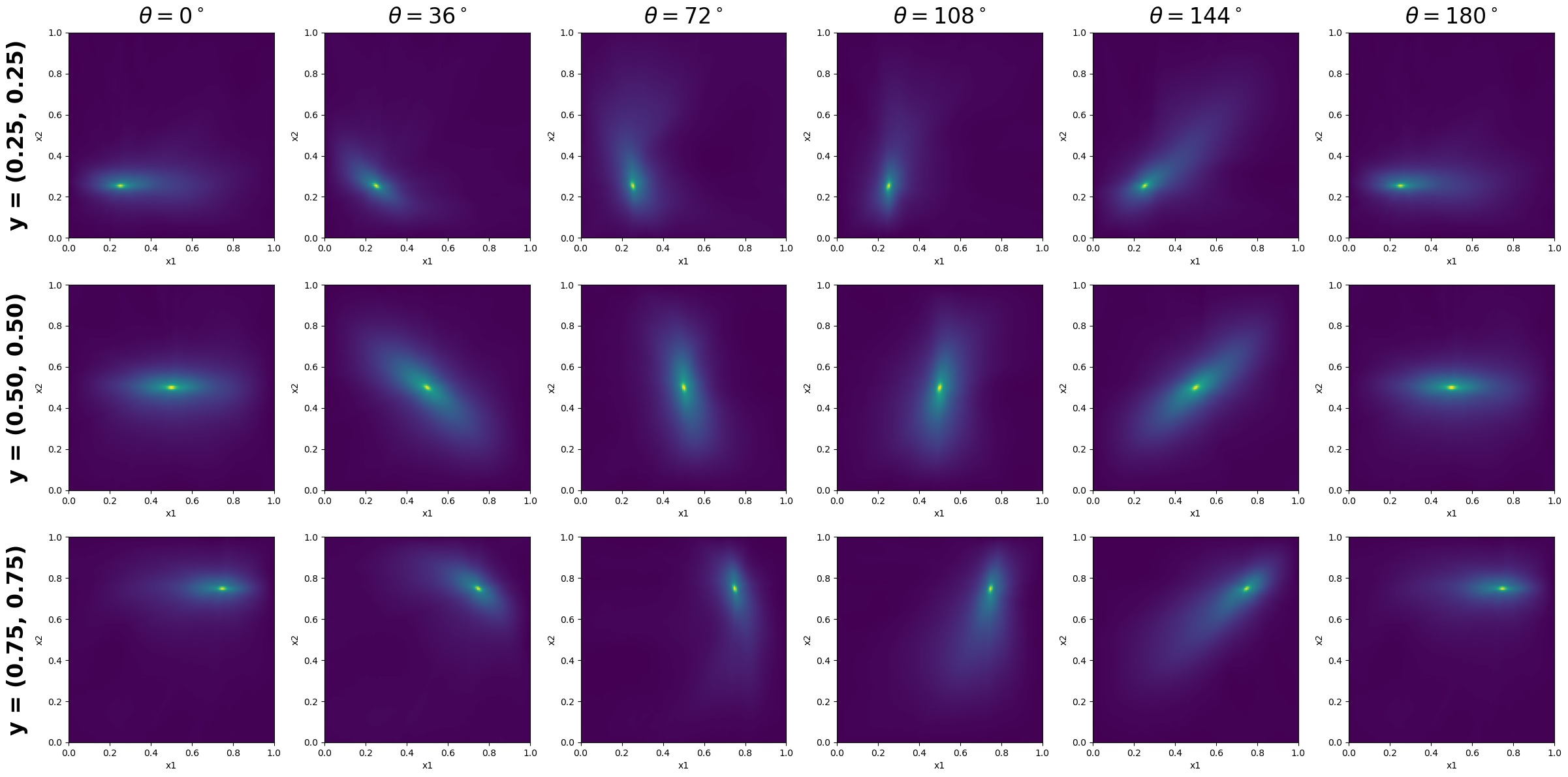}
    \caption{Parametric $\alpha$MSNN approximation to the Green's function of the rotated Laplacian operator with $\xi=0.1$ evaluated at several sources $y$ and rotation angles $\theta$.}
    \label{fig:rlap greens}
\end{figure}

Figure~\ref{fig:rlap greens} shows the learned parametric $\alpha$MSNN 
evaluated at several source locations $y$ (rows) and rotation angles $\theta$ (columns) 
for $\xi = 0.1$. 
As expected, the peak of the approximate Green’s function consistently occurs near $x \approx y$ 
and shifts appropriately with the source position. 
More importantly, as $\theta$ varies, the anisotropic lobe rotates smoothly, 
aligning with the principal direction of the diffusion tensor. 
This smooth angular interpolation is noteworthy: 
the network produces physically consistent kernels even for angles not included 
in the four anchor values used during training. 
This demonstrates that the single model has not merely memorized discrete cases 
but has learned to generalize across the continuous parameter space, 
yielding a coherent and physically faithful representation of the entire operator family.

\begin{figure}[htbp]
    \centering
    \includegraphics[width=0.95\linewidth]{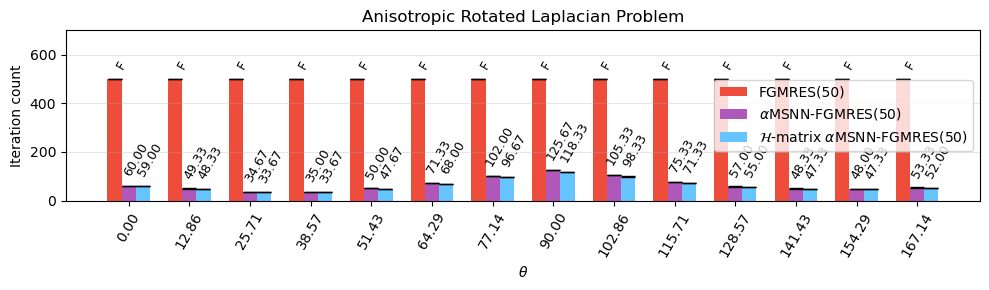}
    \caption{$\theta$ vs number of FGMRES(50) iterations for anisotropic rotated Laplacian problem with $\xi=0.1$. We discretize all problems on a $98\times98$ grid. F indicates that the solver fails to converge within $500$ iterations.}
    \label{fig:rlap iterations}
\end{figure}

To provide quantitative validation of the visually accurate kernels shown in Figure~\ref{fig:rlap greens}, we evaluate solver performance using a finer $98\times98$ grid and a more comprehensive set of rotation angles~$\theta$, \rev{including angles not used to construct the anchor data}.
Figure~\ref{fig:rlap iterations} reports the number of FGMRES(50) iterations 
required for convergence on a $98\times98$ grid. 
The unpreconditioned solver fails to converge for all tested angles, 
indicating that the problem remains challenging throughout the parameter space. 
In contrast, both dense matrix and $\mathcal{H}$-matrix preconditioners 
achieve robust and nearly uniform performance across~$\theta$. 
These results demonstrate the efficiency and generalization capability of the proposed approach: 
a single trained $\alpha$MSNN model provides an effective parametric preconditioner 
\rev{across the tested rotation angles}.
This success further underscores the potential of our framework 
for addressing broader classes of complex parameterized PDEs.

\subsection[Experiment 4: Schur Complement Preconditioning]{\rev{Experiment 4: Schur Complement Preconditioning}}
\label{subsec:schur_experiments}
{

We next evaluate the learned Green kernels within the two-level ILU framework described in Section~\ref{sec:schur_deployment}.
We consider two parametrized problems: the rotated diffusion
problem~\eqref{eq:rlap}, with $\theta\in[0,\pi)$ and
$\xi\in[0.02,1]$, and an indefinite reaction--diffusion problem
with $a(x)=I$, $\mathbf{b}(x)=0$, and $c(x)=-2\pi^2r$,
where $r\in[1.25,2]$.
Both problems use homogeneous Dirichlet conditions on the unit square.
Training concentrates on the Green kernel needed for the inverse Schur complement: source coordinates are sampled within strips of half-width $3/64$ around the central separator, while target and boundary samples cover the physical domain.

For each problem, we train a three-scale $\alpha$MSNN with
$\varepsilon\in\{0.01,0.005,0.002\}$ and three hidden layers per branch.
The far-field network has $100$ neurons per hidden layer,
while the near- and middle-field branches use $20$, $30$,
and $40$ neurons per hidden layer at the successive scales.
The finest scale uses nine experts; the scaling and gating networks have width $10$, giving $99{,}645$ parameters in total.
We use four stages of $20{,}000$ Adam updates, with $4{,}000$ residual, $1{,}000$ boundary, and $1{,}000$ anchor samples per update.
The loss weights are $w_{\rm res}=\varepsilon^2$ and $w_{\rm bc}=w_{\rm aux}=100$.
The anchors use all $32\times32$ interior target nodes and
$240$ source points near $\Gamma$, with $20$ parameter settings for rotated diffusion and five for the indefinite problem.
GPU training took $67.41$ and $66.42$ minutes for the main Adam phases, respectively.
Each trained kernel is reused on grids with $n=63,127,255$ interior points per direction, giving $n^2$ unknowns and $n_\Gamma=2n-1$ separator nodes.

We implement the two-level ILU preconditioner of Section~\ref{sec:schur_deployment} using the parGeMSLR package~\cite{Xu_Kalantzis_Li_Xi_Dillon_Saad_2022} and denote this configuration by Green--Schur. We compare it with a classical two-level ILUT configuration in the same package and a two-level Schwarz--ILUT method implemented using PETSc and its SuperLU backend.
All three configurations use four physical subdomains arranged in a $2\times2$ layout.
Green--Schur and the classical two-level ILUT configuration explicitly separate the central interface and share the same interior ILUT solves, with low-rank corrections disabled.
Each applies its separator preconditioner once per outer iteration.
The interior ILUT drop tolerance is $10^{-6}$; the classical configuration uses ILUT of $\mathbf C$ with drop tolerance $10^{-2}$, while Green--Schur uses the learned matrix $\mathbf T$.
The Schwarz method uses one-layer overlapping subdomains with local ILUT solves and a coarse grid correction on a $32\times32$ grid.
Its local ILUT drop tolerance is also $10^{-6}$, and the local fill ratios are approximately matched to those of the interior factors used by Green--Schur and two-level ILUT.

All methods use the same matrices and ten normalized right-hand sides generated from five uniform and five Gaussian random vectors.
The outer solver is the right-preconditioned FGMRES(50) implementation in parGeMSLR, with a zero initial guess, a maximum of $500$ iterations, and a relative residual tolerance of $10^{-8}$.
CPU setup and solves use one MPI rank and one thread.
Table~\ref{tab:schur_comparison} reports outer FGMRES iteration
counts (Its.) and setup and solve times in seconds for two
representative parameter settings.
The two classical preconditioners are constructed on the CPU, while Green--Schur uses the CPU for local ILUT factorization and the GPU to form $\mathbf T$.
At $n=255$, CPU setup takes $0.4220$ and $1.2472$ seconds for the rotated and indefinite problems, respectively, while warm GPU generation takes $0.5623$ and $0.5747$ seconds.
These timings suggest that most of the Green-generation cost on the largest grid can be hidden by overlapping GPU generation with CPU setup.

\begin{table}[htbp]

\centering
\caption{\rev{Comparison of two-level preconditioners. Its. and Solve are medians over ten right-hand sides; Setup reports median construction times. All times are in seconds. For Green--Schur, Setup$^*$ is the idealized estimate $\max\{t_G,t_{\rm CPU}\}$ for full CPU--GPU overlap, using separately measured warm Green generation and CPU setup times.}}
\label{tab:schur_comparison}
\small
\setlength{\tabcolsep}{3pt}
\renewcommand{\arraystretch}{1.12}
\begin{tabular*}{\linewidth}{@{\extracolsep{\fill}}r*{9}{r}@{}}
\toprule
& \multicolumn{3}{c}{Two-level ILUT}
& \multicolumn{3}{c}{Schwarz--ILUT}
& \multicolumn{3}{c}{Green--Schur} \\
\cmidrule(lr){2-4}\cmidrule(lr){5-7}\cmidrule(l){8-10}
$n$ & Its. & Setup & Solve & Its. & Setup & Solve & Its. & Setup$^*$ & Solve \\
\midrule
\multicolumn{10}{l}{Rotated diffusion: $\theta=\pi/3$, $\xi=0.03$} \\
63 & 43 & 0.0084 & 0.0117 & 12 & 0.0136 & 0.0031 & 30.5 & 0.0347 & 0.0080 \\
127 & 67 & 0.0565 & 0.1015 & 23 & 0.0599 & 0.0315 & 35 & 0.1335 & 0.0545 \\
255 & 99 & 0.4170 & 1.8611 & 69 & 0.4176 & 0.9475 & 28 & 0.5623 & 0.5175 \\
\midrule
\multicolumn{10}{l}{Indefinite: $r\approx1.4482$} \\
63 & 28 & 0.0104 & 0.0080 & 10 & 0.0131 & 0.0026 & 13 & 0.0335 & 0.0036 \\
127 & 41 & 0.1214 & 0.1096 & 18 & 0.0833 & 0.0373 & 11 & 0.1450 & 0.0276 \\
255 & 71 & 1.2517 & 2.0518 & 37 & 0.7493 & 0.6482 & 12 & 1.2472 & 0.3371 \\
\bottomrule
\end{tabular*}
\end{table}

All runs reach the prescribed tolerance.
Green--Schur requires fewer outer iterations and less solve time than classical two-level ILUT in all six cases.
At $n=255$, the rotated and indefinite solves take $0.517$ and $0.337$ seconds, compared with $1.861$ and $2.052$ seconds for two-level ILUT and $0.947$ and $0.648$ seconds for Schwarz--ILUT.
Schwarz--ILUT has a shorter setup time, but these solve-time savings can offset the higher setup cost when solving multiple right-hand sides with the same matrix.
With grid refinement, Green--Schur iteration counts level off at about $30$ for rotated diffusion and remain between $11$ and $13$ for the indefinite problem, whereas those of both classical preconditioners continue to increase.
Future work will further optimize the implementation and explore its parallel performance on distributed CPU--GPU systems.

}

\section{Conclusion}\label{sec:conclusion}
We have presented a data-driven framework for constructing efficient
approximate inverse preconditioners for diffusion-convection-reaction PDEs
by learning the Green's function using an Adaptive Multiscale Neural Network
($\alpha$MSNN) approach.
The method combines a multiscale neural representation with optional
coarse-grid anchor data, a multi-$\varepsilon$ staged training protocol,
and an overlapping domain decomposition approach
to improve stability, adaptability, and physical consistency.
Once trained, the neural Green's function can be discretized and compressed
into an $\mathcal{H}$-matrix representation,
with nearly linear storage and matrix-vector multiplication costs for bounded
block ranks, while maintaining favorable spectral properties.
\rev{It can also approximate the inverse of the Schur complement in a
two-level ILU preconditioner, avoiding construction of the full dense inverse.}
Numerical results on convection-dominated, indefinite, and anisotropic problems
demonstrate that both dense and $\mathcal{H}$-matrix preconditioners
deliver rapid GMRES convergence and strong scalability.
\rev{The two-level formulation also achieves competitive online solve times
compared with classical two-level ILUT and Schwarz preconditioners.}

\rev{Future work will extend this framework to higher-dimensional and time-dependent PDEs and investigate parallel and mixed-precision implementations for large-scale scientific computing applications.}

\section*{\texorpdfstring{\rev{Acknowledgments}}{Acknowledgments}}
\rev{The authors used GPT-5 and GPT-6 for initial literature searches, coding assistance, and refinement of the manuscript. The authors reviewed and verified the resulting text, code, and references. The authors assume responsibility for all content.}

\end{document}